\documentclass[journal]{IEEEtran}

\usepackage{graphicx}
\usepackage{graphics}
\usepackage{times}
\usepackage{amsmath}
\usepackage{amssymb}
\usepackage[ruled,vlined]{algorithm2e}
\usepackage{pseudocode}

\usepackage{enumitem}
\usepackage{cite}
\usepackage{url}

\usepackage[usenames, dvipsnames]{color}
\definecolor{darkblue}{rgb}{0.0, 0.0, 0.7}
\usepackage[colorlinks	= true,
			raiselinks	= true,
			linkcolor	= black, 
			citecolor	= blue,
			urlcolor	= ForestGreen,
			pdftitle	= {},
			pdfkeywords	= {},
			pdfsubject	= {},
			plainpages	= false]{hyperref}

\newtheorem{thm}{\bf Theorem}

\newtheorem{problem}{\bf Problem}
\newtheorem{remark}{\bf Remark}
\newtheorem{definition}{\bf Definition}
\newtheorem{claim}{\bf Claim}

\newtheorem{prop}{\bf Proposition}
\newtheorem{lem}{\bf Lemma}
\newtheorem{assumption}{\bf Assumption}

\def\QED{~\rule[-1pt]{5pt}{5pt}\par\medskip}
\newenvironment{pf}{{\bf Proof: \ }}{ \hfill \QED}

\newcommand{\calx}{\mathcal{X}}

\newcommand{\calu}{\mathcal{U}}

\DeclareMathOperator*{\argmin}{arg\,min}

\newcommand{\upd}{\color{black}}
\newcommand{\updnew}{\color{black}}
\newcommand*{\IFIEEE}{}%

\newcommand*{\LONGVERSION}{}

\title{
Transfer-Entropy-Regularized \\  Markov Decision Processes 
}

\author{Takashi Tanaka$^{1}$ \hspace{10ex}
 Henrik Sandberg$^{2} $ \hspace{10ex} Mikael Skoglund$^{3}$ 
\thanks{$^{1}$University of Texas at Austin, USA, {\tt ttanaka@utexas.edu};  $^{2}$KTH Royal Institute of Technology, Sweden, {\tt hsan@kth.se}; $^{3}$KTH Royal Institute of Technology, Sweden, {\tt skoglund@kth.se}.}%
        }

\begin{document}

\date{}
\maketitle

\begin{abstract}                          
We consider the framework of transfer-entropy-regularized Markov Decision Process (TERMDP) in which the weighted sum of the classical state-dependent cost and the transfer entropy from the state random process to the control random process is minimized.   
Although TERMDPs are generally formulated as nonconvex optimization problems, we derive an analytical necessary optimality condition expressed as a finite set of nonlinear equations, based on which an iterative forward-backward computational procedure similar to the Arimoto-Blahut algorithm is proposed. 
{\upd It is shown that every limit point of the sequence generated by the proposed algorithm is a stationary point of the TERMDP.}
Applications of TERMDPs are discussed in the context of networked control systems theory and non-equilibrium thermodynamics.
The proposed algorithm is applied to an information-constrained maze navigation problem, whereby we study how the price of information qualitatively alters the optimal decision polices.
\end{abstract}

\section{Introduction}
{\updnew
\emph{Transfer entropy} \cite{schreiber2000measuring} is a quantity that can be understood as a measure of information flow between random processes. 
It is a generalization of \emph{directed information}, a concept proposed in the information theory literature for the analysis of communication systems with feedback \cite{marko1973bidirectional,massey1990causality,kramer1998directed}. 
Closely related concepts include the \emph{KL-causality measure}  \cite{gourieroux1987kullback}, which was originally introduced in the economic statistics literature for the causality analysis.\footnote{Sometimes (e.g., in statistical physics \cite{parrondo2015thermodynamics}), transfer entropy is used as a synonym for directed information. It appears that the concepts of transfer entropy \cite{schreiber2000measuring}, directed information \cite{marko1973bidirectional,massey1990causality}, and Kullback causality measure \cite{gourieroux1987kullback} were introduced independently.} Recently, these concepts have been applied in a broad range of academic disciplines, including neuroscience \cite{wibral2014directed}, finance \cite{jiao2013universal}, and social science \cite{ver2012information}.

In this paper, we formulate the problem of \emph{transfer-entropy-regularized Markov Decision Process (TERMDP)}, and develop a numerical solution algorithm. TERMDP is an optimal control problem in which we seek a causal decision-making policy that minimizes the weighted sum of the classical state-dependent cost and transfer entropy from the state random process to the control actions.
As we will discuss in the sequel, TERMDP predicts a fundamental performance limitation of feedback control systems from an information-theoretic perspective. 
The first context in which TERMDP naturally arises is networked control systems theory, where the trade-off between the best achievable control performance and the data rate at which sensor information is fed back to the controller is a central question.
Prior work has shown that transfer entropy can be used as a proxy for the data rate on communication channels, and thus solving TERMDP provides a fundamental performance limitation of such systems. 
The second application of TERMDP is non-equilibrium thermodynamics. There has been renewed interests in the generalized second law of thermodynamics, in which transfer entropy arises as a key concept \cite{ito2013information}.
TERMDP in this context can be interpreted as the problem of operating thermal engines at a nonzero work rate near the fundamental limitation of the second law of thermodynamics.

In contrast to the standard MDP  \cite{puterman2014markov}, TERMDP penalizes the  information flow from the underlying state random process to the control random process. Consequently, TERMDP promotes  ``information-frugal'' decision policies, under which control actions tend to be statistically less dependent on the underlying Markovian state dynamics. 
This is often a favorable property in various real-time decision-making scenarios (for both humans and robots) in which information acquisition, processing, and transmission are costly operations. 
Therefore, it is expected that TERMDP plays major roles in broader contexts beyond the aforementioned applications, 
although the interpretations of transfer entropy in each application must be carefully discussed. 

In the literature, a few alternative approaches have been suggested to apply information-theoretic cost functions to capture decision-making costs in MDPs. Similarities and differences between TERMDP and the existing problem formulations are noteworthy.
The \emph{rationally inattentive control} problem \cite{sims2003implications, shafieepoorfard2016rationally} has been motivated in a macroeconomic context, where Shannon's mutual information (a special case of transfer entropy) is adopted as an attention cost for decision-makers.
The authors of \cite{todorov2007linearly,theodorou2010generalized,dvijotham2012unifying} present a class of optimal control problems in which control costs are modeled as the Kullback-Leibler (KL) divergence from the ``uncontrolled'' state trajectories to the ``control'' state trajectories. 
Alternative information-theoretic decision costs in dynamic environments include \emph{predictive information} \cite{bialek2001predictability}, \emph{past-future information-bottleneck} \cite{creutzig2009past}, and \emph{information-to-go}  \cite{tishby2010information, rubin2012trading}.
Information-theoretic bounded rationality and its analogy to thermodynamics are discussed in \cite{ortega2013thermodynamics}. 
While intuitively plausible, some of these problem formulations lack physical (or coding-theoretic) justifications, unlike TERMDP, whose operational interpretation can be found in the aforementioned contexts. 

An equivalent problem formulation to TERMDP first appeared in \cite{charalambous2014optimization} and \cite{stavrou2015information}, where the problem was formulated in a general (Polish state space) setup.
Linear-Quadratic-Gaussian (LQG) control with minimum directed information \cite{tanaka2015lqg} is a version of TERMDP specialized to the LQG regime.
While the problem in the LQG setup was shown to be tractable by semidefinite programming \cite{tanaka2015lqg}, algorithmic aspects of TERMDP beyond the LQG regime have not been thoroughly studied.
Therefore, the primary goal of this paper is to provide an efficient computational algorithm to find a stationary point (an optimal solution candidate) of the given TERMDP. The contributions of this paper are as follows:
\begin{itemize}
\item We derive a necessary optimality condition expressed as a set of nonlinear equations involving a finite number of variables. This result recovers, and partly strengthens, results obtained in prior work\cite{charalambous2014optimization}.
\item We propose a forward-backward iterative algorithm that can be viewed as a generalization of the Arimoto-Blahut algorithm \cite{blahut1972computation,arimoto1972algorithm} to solve the optimality condition numerically. 
\end{itemize}
The proposed algorithm is the first application of the the Arimoto-Blahut algorithm for transfer entropy \emph{minimization}. Our algorithm should be compared with the generalized Arimoto-Blahut algorithm for transfer entropy \emph{maximization} proposed in \cite{naiss2013extension}. The algorithm in \cite{naiss2013extension} can be viewed as a generalization of the Arimoto-Blahut ``capacity algorithm'' in \cite{blahut1972computation}, while our proposed algorithm can be viewed as a generalization of the Arimoto-Blahut ``rate-distortion algorithm'' in   \cite{blahut1972computation}.
Unfortunately, we discover that the proposed algorithm may not converge to the global minimum due to the non-convex nature of TERMDP. This result is somewhat surprising as the global convergence of the original Arimoto-Blahut rate-distortion algorithm, which is a special case of our algorithm, is well-known.
Nevertheless, observing that the proposed algorithm belongs to the class of block coordinate descent (BCD) algorithms, we show that every limit point generated by the algorithm is guaranteed to be a stationary point of the given TERMDP.
}

Organization of the paper:
The problem formulation of TERMDP is formally introduced 
in Section~\ref{secformulation}. Mathematical preliminaries are summarized in Section~\ref{secprelim}.
Section~\ref{secmain} presents the main results. Derivation of the main results are summarized in Section~\ref{secderivation}.
Section~\ref{secinterpret} discusses applications of the TERMDP framework. A numerical demonstration of the proposed algorithm is presented in Section~\ref{secnum}. We conclude with a list of future work in Section~\ref{secsummary}.

Notation: Upper case symbols such as $X$ are used to represent random variables, while lower case symbols such as $x$ are used to represent a specific realization. Notation $x_k^l\triangleq (x_k, x_{k+1}, ... , x_l)$ and $x^t\triangleq (x_1, x_2, ... , x_t)$ will be used to specify subsequences.
 We use the natural logarithm $\log(\cdot)=\log_e (\cdot)$ throughout the paper.

\section{Problem formulation}
\label{secformulation}
We formulate TERMDP based upon
the standard Markov Decision Process (MDP) formalism  \cite{puterman2014markov} defined by a time index $t=1, 2, ... , T$, state space $\mathcal{X}_t$,  action space $\mathcal{U}_t$, transition probability $p_{t+1}(x_{t+1}|x_t, u_t)$, cost functions  $c_t: \mathcal{X}_t \times \mathcal{U}_t \rightarrow \mathbb{R}$ for each $t=1, 2, ... , T$ and $c_{T+1}: \mathcal{X}_{T+1}  \rightarrow \mathbb{R}$.
For simplicity, we assume that both $\mathcal{X}_t$ and  $\mathcal{U}_t$ are finite.
The decision policy to be synthesized can be probabilistic and history-dependent in general, and is  represented by a conditional probability distribution:
\begin{equation}
\label{eqpolicyinf}
q_t(u_t|x^t, u^{t-1}).
\end{equation}
The joint distribution of the state and control trajectories is denoted by $\mu_{t+1}(x^{t+1},u^t)$, which is 
uniquely determined by the initial state distribution $\mu_1(x_1)$, the state transition probability $p_{t+1}(x_{t+1}|x_t, u_t)$ and the decision policy $q_t(u_t|x^t, u^{t-1})$ by a recursive formula
\ifdefined\IFONECOL
\begin{equation}
\mu_{t+1}(x^{t+1}, u^{t})=p_{t+1}(x_{t+1}|x_t,u_t)q_t(u_t|x^t,u^{t-1})\mu_t(x^t,u^{t-1}). \label{eqmupq}
\end{equation}
\else
\begin{align}
&\mu_{t+1}(x^{t+1}, u^{t}) \nonumber \\
&=p_{t+1}(x_{t+1}|x_t,u_t)q_t(u_t|x^t,u^{t-1})\mu_t(x^t,u^{t-1}). \label{eqmupq}
\end{align}
\fi
{\updnew A stage-additive cost functional
\begin{equation}
\label{equsualmdp}
J(X^{T+1}, U^T) \triangleq  \sum_{t=1}^T \mathbb{E} c_t(X_t, U_t)+\mathbb{E}  c_{T+1}(X_{T+1})
\end{equation}
 is a function of random variables $X^{T+1}$ and $U^T$ with a joint distribution $\mu_{T+1}(x^{T+1},u^T)$.}
Transfer entropy is an information-theoretic quantity defined as follows:
\begin{definition}
For nonnegative integers $m$ and $n$, 
the \emph{transfer entropy of degree} $(m,n)$ is defined by
{\upd
\begin{align}
& I_{m,n}(X^T \rightarrow U^T) \triangleq \sum_{t=1}^T \mathbb{E}\log\frac{\mu_{t+1}(U_t|X_{t-m}^t, U_{t-n}^{t-1})}{\mu_{t+1}(U_t|U_{t-n}^{t-1})} \label{eqdefte} \\
&=\sum_{t=1}^T \sum_{x^t \in \mathcal{X}^t, u^t \in \mathcal{U}^t} \mu_{t+1}(x^t, u^t)\log\frac{\mu_{t+1}(u_t|x_{t-m}^t, u_{t-n}^{t-1})}{\mu_{t+1}(u_t|u_{t-n}^{t-1})}.
\nonumber
 \end{align}
 }
\end{definition}
Using conditional mutual information  \cite{CoverThomas},  transfer entropy can also be written as
\begin{equation}
 I_{m,n}(X^T \rightarrow U^T) \triangleq \sum_{t=1}^T I(X_{t-m}^t;U_t|U_{t-n}^{t-1}).
 \end{equation}
When $m=\infty$ and $n=\infty$, \eqref{eqdefte} coincides with \emph{directed information} \cite{massey1990causality}:
\begin{equation}
 I(X^T \rightarrow U^T) \triangleq \sum_{t=1}^T I(X^t;U_t|U^{t-1}). \label{eqdefdi}
 \end{equation}
The main problem studied in this paper is now formulated as follows.

\begin{problem} (TERMDP)
Let the initial state distribution $\mu_1(x_1)$ and the state transition probability $p_{t+1}(x_{t+1}|x_t, u_t)$ be given, and assume that the joint distribution $\mu_{t+1}(x^{t+1}, u^t)$ is recursively given by \eqref{eqmupq}. 
For a fixed constant $\beta\geq 0$, the \emph{Transfer-Entropy-Regularized Markov Decision Processes} is the optimization problem 
\begin{equation}
\min_{\{q_t(u_t|x^t, u^{t-1})\}_{t=1}^T}  J(X^{T+1}, U^T) + \beta  I_{m,n}(X^T \rightarrow U^T).  \label{eqmainprob}
\end{equation}
\end{problem}

A few remarks are in order regarding this problem formulation.
First, the transfer entropy term in \eqref{eqmainprob} is  interpreted as an additional cost corresponding to the information transfer from the state random process $X_t$ to the control random process $U_t$. 
The regularization parameter $\beta\geq 0$ can be thought of as the cost of information transfer.  When $\beta = 0$, the standard MDP formulation is recovered. As we increase $\beta>0$, the optimal decision policy for \eqref{eqmainprob} tends to be  ``information frugal'' in order to reduce  $ I_{m,n}(X^T \rightarrow U^T)$. That is, control actions generated by the policy becomes statistically less dependent on the state of the system.

Second, in contrast to the standard MDP (i.e., $\beta=0$)
for which optimal polities can be assumed deterministic and history-independent (Markovian) without loss of generality \cite[Section 4.4]{puterman2014markov}, {\updnew TERMDP \eqref{eqmainprob} with $\beta>0$ admits an optimal policy that is randomized and history-dependent.}
Thus, the cardinality of the solution space we must explore to solve \eqref{eqmainprob} is much larger than that of the standard MDP.
However, in Proposition~\ref{theooptcond} below, we show that one can assume without loss of performance a structure of the optimal policy of the form
\begin{equation}
\label{eqpolicy}
q_t(u_t|x_t, u^{t-1}_{t-n})
\end{equation}
instead of \eqref{eqpolicyinf}. In other words, it is sufficient to consider a policy that is  dependent only on the most recent realization of the state and the last $n$ realizations of the control inputs.

Finally, the structure of the problem \eqref{eqmainprob} is similar {\updnew (but not equivalent)} to that of the KL control formulation in \cite{todorov2007linearly}. In particular, if $(m,n)=(0,0)$, the transfer entropy cost \eqref{eqdefte} becomes 
 \[
{\upd\sum\nolimits_{t=1}^T \mathbb{E}\log\frac{\mu_{t+1}(U_t|X_t)}{\mu_{t+1}(U_t)}.}
\]
On the other hand, the KL control framework considers the KL divergence cost of the form
\[
{\upd
\sum\nolimits_{t=1}^T \mathbb{E}\log\frac{\mu_{t+1}(U_t|X_t)}{r_{t+1}(U_t|X_t)} }
\]
where  $r_{t+1}(u_t|x_t)$ is the conditional distribution specified by a predefined ``reference'' policy. 
{\updnew 
Unfortunately, this difference renders \eqref{eqmainprob} nonconvex as we will observe in Section~\ref{secnonconvex}.
Despite this disadvantage, we emphasize the importance of studying TERMDP as it arises in some practical problems in science and engineering as discussed in Section~\ref{secinterpret}. }

\section{Preliminaries}
\label{secprelim}
In this section, we summarize some technical results needed to derive our main results in this paper. 
\subsection{Structure of the optimal solution}
We first derive some important structural properties of the optimal policy, which allow us to rewrite the main problem \eqref{eqmainprob} in a simpler form.
The desired structural results can be obtained by applying the dynamic programming principle to  \eqref{eqmainprob}.
To this end, notice that \eqref{eqmainprob} can be viewed as a $T$-stage optimal control problem, in which the joint distribution $\mu_t$ is the ``state'' of the system to be controlled by a multiplicative control action $q_t$ via the state evolution equation \eqref{eqmupq}. Introduce the value function as
\begin{align}
&V_k\left(\mu_k(x^k, u^{k-1})\right) \triangleq  \nonumber\\
&\min_{\{q_t\}_{t=k}^T}   \sum_{t=k}^T \left\{ \mathbb{E} {\upd c_t(X_t, U_t)}+I(X_{t-m}^t;U_t|U_{t-n}^{t-1})\right\}. \label{eqvk}
\end{align}
The value function satisfies the Bellman equation
\begin{align}
&V_t\left(\mu_t(x^t, u^{t-1})\right) 
=\min_{q_t} \Bigl\{ \mathbb{E}c_t(X_t, U_t) \biggr.  \nonumber \\
& \hspace{5ex}\Bigl. + I(X_{t-m}^t;U_t|U_{t-n}^{t-1})  + V_{t+1}(\mu_{t+1}(x^{t+1}, u^{t}))\Bigr\} \label{eqhjb1}
\end{align}
for $t=1,2, ... , T$, with the terminal condition 
\begin{equation}
V_{T+1}\left(\mu_{T+1}(x^{T+1}, u^{T})\right) =\mathbb{E}^{\mu_{T+1}} c_{T+1}(X_{T+1}). \label{eqterminal}
\end{equation}
The next proposition summarizes key structural results.
\begin{prop}
\label{propstructure}
For the optimization problem \eqref{eqmainprob} and its dynamic programming formulation \eqref{eqvk}--\eqref{eqterminal}, the following statements hold for each $k=1, 2, ... , T$.
\begin{itemize}
\item[(a)] For each sequence of policies $\{q_t(u_t|x^t, u^{t-1})\}_{t=k}^T$, there exists a sequence of policies of the form $\{q'_t(u_t|x_t, u_{t-n}^{t-1})\}_{t=k}^T$ such that the value of the objective function on the left hand side of  
\eqref{eqvk} attained by $\{q'_t\}_{t=k}^T$ is less than or equal to the value  attained by $\{q_t\}_{t=k}^T$.
\item[(b)] If the policy at time step $k$ is of the form $q'_k(u_k|x_k, u_{k-n}^{k-1})$, the identity
$I(X_{k-m}^k;U_k|U_{k-n}^{k-1})=I(X_k;U_k|U_{k-n}^{k-1})$
holds.
\item[(c)] The value function $V_k(\mu_k(x^k, u^{k-1}))$ depends only on the marginal distribution $\mu_k(x_k, u_{k-n}^{k-1})$.
\end{itemize}
\end{prop}
\begin{pf}
\ifdefined\LONGVERSION
See Appendix~\ref{appstructure}.
\else
See [xx, Appendix A].
\fi
\end{pf}

An implication of Proposition~\ref{propstructure}~(a) is that the search for the optimal policy for the original TERMDP \eqref{eqmainprob} can be restricted to the class of policies of the form $q_t(u_t|x_t, u_{t-n}^{t-1})$ without loss of performance. 
Proposition~\ref{propstructure} (b) implies that, as far as the policy of the form $q_t(u_t|x_t, u_{t-n}^{t-1})$ is used,  the problem remains equivalent even after the transfer entropy term $I_{m,n}(X^T\rightarrow U^T)$ is replaced by the transfer entropy of degree $(0, n)$:
\[
I_{0,n}(X^T\rightarrow U^T)=\sum\nolimits_{t=1}^T I(X_t;U_t|U_{t-n}^{t-1}).
\]
Proposition~\ref{propstructure}~(c) implies that the distribution $\mu_t(x_t, u_{t-n}^{t-1})$, rather than the original  $\mu_t(x^t, u^{t-1})$, suffices as the state of the considered problem. This ``reduced'' state evolves according to
\begin{align}
&\mu_{t+1}(x_{t+1}, u_{t-n+1}^t)=\nonumber \\
&\sum_{x_t \in \mathcal{X}_t, u_{t-n}\in \mathcal{U}_{t-n}}p_{t+1}(x_{t+1}|x_t, u_t)q_t(u_t|x_t, u_{t-n}^{t-1})\mu_t(x_t, u_{t-n}^{t-1}). \label{eqmustate}
\end{align}
Based on these observations, it can be seen that Problem 1 can be solved by solving the following simplified problem:
\begin{problem} (Simplified TERMDP)
Let the initial state distribution $\mu_1(x_1)$ and the state transition probability $p_{t+1}(x_{t+1}|x_t, u_t)$ be given, and assume that the joint distribution $\mu_t(x_t, u_{t-n}^{t-1})$ is recursively given by \eqref{eqmustate}. 
For a fixed constant $\beta\geq 0$, the \emph{simplified TERMDP} is the optimization problem 
\begin{equation}
\min_{\{q_t(u_t|x_t, u_{t-n}^{t-1})\}_{t=1}^T}  J(X^{T+1}, U^T) + \beta  I_{0,n}(X^T \rightarrow U^T).  \label{eqsimpleprob}
\end{equation}
\end{problem}

{\updnew Notice that Proposition~\ref{propstructure} implies that a global minimizer for \eqref{eqsimpleprob}  is a global minimizer for \eqref{eqmainprob}. With an appropriate notion of local optimality, it can also be shown that a local minimizer for \eqref{eqsimpleprob} is also a local minimizer for \eqref{eqmainprob}, as detailed in 
\ifdefined\LONGVERSION
Appendix~\ref{apdxlocal}.
\else
[xx, Appendix B].
\fi
}
For this reason, in what follows, we will develop an algorithm that solves  the simplified TERMDP \eqref{eqsimpleprob} rather than the original TERMDP \eqref{eqmainprob}.

Proposition~\ref{propstructure} implies that an optimal solution to both the original and simplified TERMDP can be found by solving the Bellman equation
\begin{align}
&V_t\left(\mu_t(x_t, u_{t-n}^{t-1})\right) 
=\min_{q_t} \Bigl\{ \mathbb{E}c_t(X_t, U_t) \biggr.  \nonumber \\
& \hspace{5ex}\Bigl. + I(X_t;U_t|U_{t-n}^{t-1})  + V_{t+1}(\mu_{t+1}(x_{t+1}, u_{t-n+1}^{t}))\Bigr\} \label{eqhjb2}
\end{align}
with the state transition rule \eqref{eqmustate} and the terminal condition
\begin{equation}
V_{T+1}\left(\mu_{T+1}(x_{T+1}, u_{T-n+1}^{T})\right) =\mathbb{E}^{\mu_{T+1}} c_{T+1}(X_{T+1}). \label{eqterminal2}
\end{equation}
Notice that the Bellman equation \eqref{eqhjb2} is simpler than the original form \eqref{eqhjb1}.
However, solving \eqref{eqhjb2} remains computationally challenging as the right hand side of \eqref{eqhjb2} involves a nonconvex optimization problem. In Section~\ref{secnonconvex}, we present a simple numerical example demonstrating this nonconvexity.

\subsection{Transfer entropy and directed information}
\label{secbound}
In some applications (e.g., networked control systems, see Section~\ref{secncs}), we are interested in \eqref{eqmainprob} with directed information ($m=\infty$ and $n=\infty$), even though solving such a problem is often computationally intractable. In such applications, approximating directed information with transfer entropy with finite degrees is a natural idea.
{\updnew The next proposition describes the consequence of such an approximation.}
\begin{prop}
\label{propfinitete}
For any fixed decision policy of the form $q_t(u_t|x_t, u_{t-n}^{t-1})$, $t=1, 2, ... , T$, we have\footnote{$I_{m,n}$ is a short-hand notation for  $I_{m,n}(X^T \rightarrow U^T)$.}
\[
I_{0,n}\geq 
I_{0,n+1}\geq \cdots \geq
I_{0,\infty}.
\]
\end{prop}
\begin{pf}
\ifdefined\LONGVERSION
See Appendix~\ref{apptemonotone}.
\else
See [xx, Appendix C].
\fi
\end{pf}
The following chain of inequalities shows that the optimal value of \eqref{eqsimpleprob} with any finite $n$ provides an upper bound on the optimal value of \eqref{eqmainprob} with $(m,n)=(\infty, \infty)$.
\begin{subequations}

\begin{align}
& \min_{\{q_t(u_t|x^t, u^{t-1})\}_{t=1}^T} J(X^{T+1}, U^T) + \beta  I_{\infty, \infty} \nonumber \\
&= \min_{\{q_t(u_t|x_t, u^{t-1})\}_{t=1}^T} J(X^{T+1}, U^T)+ \beta  I_{\infty, \infty} \label{eqdivste1} \\
&= \min_{\{q_t(u_t|x_t, u^{t-1})\}_{t=1}^T} J(X^{T+1}, U^T) + \beta  I_{0,\infty} \label{eqdivste2} \\
&\leq \min_{\{q_t(u_t|x_t, u_{t-n}^{t-1})\}_{t=1}^T} J(X^{T+1}, U^T)+ \beta  I_{0,\infty} \label{eqdivste3} \\
&\leq \min_{\{q_t(u_t|x_t, u_{t-n}^{t-1})\}_{t=1}^T} J(X^{T+1}, U^T) + \beta  I_{0,n}. \label{eqdivste4}
\end{align}
\end{subequations}
Equalities \eqref{eqdivste1} and \eqref{eqdivste2} follows from Proposition~\ref{propstructure} (a) and (b), respectively. The inequality \eqref{eqdivste3} is trivial since any policy of the form $q_t(u_t|x_t, u_{t-n}^{t-1})$ is a special case of the policy of the form $q_t(u_t|x_t, u^{t-1})$.
The final inequality \eqref{eqdivste4} is due to Proposition~\ref{propfinitete}.

\subsection{Rate-distortion theory and Arimoto-Blahut Algorithm}
\label{secrd}

In the special case with $T=1$, $n=0$, $\beta=1$ and $c_{T+1}(\cdot)=0$, the optimization problem \eqref{eqsimpleprob} becomes 
\begin{equation}
\label{eqrdprob}
\min_{q(u|x)} \mathbb{E}c({\upd X,U})+I(X;U)
\end{equation}
where the probability distribution $p(x)$ on $\mathcal{X}$ is given. In this special case, the problem \eqref{eqrdprob} is convex, and the solution is well-known in the context of  rate-distortion theory \cite{CoverThomas}.
\begin{prop}
\label{propext}
A conditional distribution $q^*(u|x)$ is a {\updnew global minimizer for} \eqref{eqrdprob} if and only if it {\updnew satisfies} the  following condition $p(x)$-almost everywhere:
\begin{subequations}
 \label{eqopt}
\begin{align}
q^*(u|x)&=\frac{\nu^*(u)\exp\left\{-c(x,u)\right\}}{\sum_{u \in\calu} \nu^*(u)\exp\left\{-c(x,u)\right\}} \label{eqopt1} \\
\nu^*(u)&=\sum\nolimits_{x\in \calx} p(x)q^*(u|x). \label{eqopt2}
\end{align}
\end{subequations}
\end{prop}
\begin{pf}
This result is standard and hence the proof is omitted.
See \cite[Appendix A]{petersen2012robust} and \cite{theodorou2015nonlinear} for relevant discussions.
\end{pf}
Condition~\eqref{eqopt} is required only $p(x)$-almost everywhere since for $x$ such that $p(x)=0$, $q^*(u|x)$ can be chosen arbitrarily. Commonly, the denominator in \eqref{eqopt1} is called the \emph{partition function}:
\[
\phi^*(x)\triangleq \sum\nolimits_{u \in \calu} \nu^*(u)\exp\left\{-c(x,u)\right\}.
\]
By substitution, it is easy to show that
the optimal value of \eqref{eqrdprob} can be written in terms of $\nu^*(u)$ as
\begin{equation}
-\sum\nolimits_{x \in \calx} p(x)\log \left\{\sum\nolimits_{u \in \calu} \nu^*(u) \exp\{-c(x,u)\}\right\},
\end{equation}
or more compactly as $\mathbb{E}^{p(x)}\{-\log \phi^*(X)\}$. This quantity is often referred to as free energy \cite{parrondo2015thermodynamics,theodorou2012relative}.

The Arimoto-Blahut algorithm is an iterative algorithm to compute $q^*(u|x)$ satisfying \eqref{eqopt} numerically. It is based on the alternating updates:
\begin{subequations}
\label{eqarimotoblahut}
\begin{align}
\nu^{(k)}(u)&=\sum\nolimits_{x \in\calx} p(x)q^{(k-1)}(u|x) \\
q^{(k)}(u|x)&= \frac{\nu^{(k)}(u)\exp\{-c(x,u)\}}{\sum_{u \in \calu} \nu^{(k)}(u) \exp\{-c(x,u)\}}.
\end{align}
\end{subequations}
The algorithm is first proposed for the computation of channel capacity \cite{arimoto1972algorithm} and for the computation of rate-distortion functions \cite{blahut1972computation}. Clearly, the optimal solution $(q^*, \nu^*)$ is a fixed point of the algorithm \eqref{eqarimotoblahut}. Under a mild assumption, convergence of the algorithm is guaranteed; see  \cite{arimoto1972algorithm,csiszar1974computation, tseng2001convergence}. The main algorithm we propose in this paper to solve the simplified TERMDP  \eqref{eqsimpleprob} can be thought of as a generalization of the Arimoto-Blahut ``rate-distortion algorithm.''

\subsection{Block Coordinate Descent Algorithm}
The Arimoto-Blahut algorithm can be viewed as a block Coordinate Descent (BCD) algorithm applied to a special class of objective functions.
In this subsection, we summarize elements of the BCD method and a version of its convergence results that is relevant to our analysis. Consider the problem
\begin{subequations}
\label{eqbcdprob}
\begin{align}
\min & \hspace{2ex} f(x) \\
\text{s.t. } & \hspace{2ex} x\in X=X_1\times X_2 \times ... \times X_N
\end{align}
\end{subequations}
where the feasible set $X$ is the Cartesian product of closed, nonempty and convex subsets $X_i \subseteq \mathbb{R}^{n_i}$ for $i=1, 2, ... , N$, and the function $f: \mathbb{R}^{n_1+...+n_N}\rightarrow \mathbb{R}\cup \{\infty\}$ is  continuously differentiable on the sublevel set $\{x\in X: f(x)\leq f(x^{(0)})\}$, where $x^{(0)}\in X$ is a given initial point. We call $x^*\in X$ a \emph{stationary point} for \eqref{eqbcdprob} if it satisfies $\nabla f(x^*)^\top (y-x^*)\geq 0$ for every $y\in X$, where $\nabla f(x^*)$ is the gradient of $f$ at $x^*$. If $f$ is convex, every stationary point is a global minimizer of $f$.
The BCD algorithm for \eqref{eqbcdprob} is defined by the following cyclic update rule:
\begin{subequations}
\label{eqbcd}
\begin{align}
x_1^{(k)} &\in \argmin_{x_1} f(x_1, x_2^{(k-1)}, ... , x_N^{(k-1)}) \\
x_2^{(k)}&\in \argmin_{x_2} f(x_1^{(k)}, x_2, x_3^{(k-1)}, ... , x_N^{(k-1)}) \\
&\hspace{5ex}\cdots \nonumber \\
x_N^{(k)}&\in \argmin_{x_N} f(x_1^{(k)}, x_2^{(k)}, ... , x_N).
\end{align}
\end{subequations}
A number of sufficient conditions for the convergence of the BCD algorithm are known in the literature (e.g., \cite{bertsekas2016nonlinear, grippo2000convergence, tseng2001convergence} and references therein).
For instance, if $f$ is pseudoconvex and has compact level sets, then every limit point of the sequence $\{x^{(k)}\}$ generated by the BCD algorithm is a global minimizer of $f$ \cite[Proposition 6]{grippo2000convergence}. This result can be applied to show the global convergence of the Arimoto-Blahut algorithm \eqref{eqarimotoblahut}, simply by noticing that the objective function in \eqref{eqrdprob} can be written as a convex function of $\nu$ and $q$ as
\[
f(\nu, q)=\sum_{x\in \mathcal{X}, u\in \mathcal{U}}p(x)q(u|x)\left(c(x,u)+\log\frac{q(u|x)}{\nu(u)}\right)
\]
and that \eqref{eqarimotoblahut} is equivalent to the BCD update rule \eqref{eqbcd}.
In the absence of the convexity assumption on $f$, it is typically required that each coordinate-wise minimization  is \emph{uniquely}\footnote{The counterexample by Powell \cite{powell1973search} with $N=3$ shows that the lack of uniqueness of the coordinate-wise minimizer can result in a BCD algorithm with a limit point which is not a stationary point.} attained in order to guarantee that every limit point of the BCD algorithm is a stationary point \cite[Proposition 2.7.1]{bertsekas2016nonlinear}. 
Unfortunately, the generalized Arimoto-Blahut algorithm we introduce in this paper for the TERMDP is a BCD algorithm applied to a nonconvex objective function, and the uniqueness of the coordinate-wise minimizer cannot be assumed. Thus, none of the above results are applicable to prove the convergence.
Fortunately, the requirement of the uniqueness of the coordinate-wise minimizer can be relaxed when $N=2$ (two-block BCD algorithms). The following result is due to \cite[Corollary 2]{grippo2000convergence} and \cite[Theorem 4.2 (c)]{tseng2001convergence}.
\begin{lem}
\label{lembcd}
Consider the problem \eqref{eqbcdprob} with $N=2$, and suppose that the sequence $\{x^{(k)}\}$ generated by the two-block BCD algorithm \eqref{eqbcd} has limit points. Then, every limit point $x^*$ of $\{x^{(k)}\}$ is a stationary point of the problem \eqref{eqbcdprob}.
\end{lem}

{\updnew 
Lemma~\ref{lembcd} is critical to obtain one of our main results (Theorem~\ref{theoconvergefbaba}) below.
}

\begin{figure*}[t]
\hrule
\begin{subequations}
\label{eqoptcond}
\begin{align}
\mu^*_{t+1}(x_{t+1}, u^{t}_{t-n+1})&=\sum_{x_t \in \mathcal{X}_t}\sum_{u_{t-n}\in\mathcal{U}_{t-n}}  p_{t+1}(x_{t+1}|x_t,u_t) q^*_t(u_t|x_t, u^{t-1}_{t-n})\mu^*_t(x_t, u^{t-1}_{t-n}) \label{eqoptcond1} \\
\nu^*_t(u_t|u^{t-1}_{t-n})&=\sum_{x_t \in \mathcal{X}_t} q^*_t(u_t|x_t, u^{t-1}_{t-n})\mu^*_t(x_t|u_{t-n}^{t-1}), \;\; \forall{u_{t-n}^{t-1}} \text{ such that } \mu_t(u^{t-1}_{t-n})>0 \label{eqoptcond2} \\
\rho^*_t(x_t, u_{t-n+1}^t)&= c_t(x_t,u_t)- \sum_{x_{t+1}\in\calx_{t+1}} p_{t+1}(x_{t+1}|x_t,u_t)\log \phi^*_{t+1}(x_{t+1}, u_{t-n+1}^{t}) \label{eqoptcond3} \\
\phi^*_t(x_t, u_{t-n}^{t-1}) &= \sum_{u_t \in \mathcal{U}_t}\nu^*_t(u_t|u^{t-1}_{t-n})\exp \left\{-\rho^*_t(x_t, u_{t-n+1}^t)\right\} \label{eqoptcond4} \\
q_t^*(u_t|x_t, u_{t-n}^{t-1}) &=\frac{\nu^*_t(u_t|u^{t-1}_{t-n})\exp \left\{-\rho^*_t(x_t, u_{t-n}^t)\right\}}{\phi^*_t(x_t, u_{t-n}^{t-1})}, \;\; \forall (x_t, u_{t-n}^{t-1}) \text{ such that } \mu_t(x_t, u_{t-n}^{t-1})>0 \label{eqoptcond5}
\end{align}
\end{subequations}
\hrule
\end{figure*}

\section{Main Results}
\label{secmain}
{\updnew This section summarizes the main results (Theorems~\ref{theooptcond} and \ref{theoconvergefbaba}) of this paper. Derivations are deferred to Section~\ref{secderivation}.}

\subsection{Necessary Optimality Condition}
The first technical result states that a necessary optimality condition for the simplified TERMDP \eqref{eqsimpleprob} is given by the nonlinear condition \eqref{eqoptcond} in terms of $(\mu^*, \nu^*, \rho^*, \phi^*, q^*)$.  

\begin{thm}
\label{theooptcond}
If  $\{q_t^*\}_{t=1}^T$ is a local minimizer for \eqref{eqsimpleprob}, then there exist variables $\{\mu_{t+1}^*, \nu_t^*, \rho_t^*, \phi_t^*\}_{t=1}^T$ satisfying the set of nonlinear equations \eqref{eqoptcond} together with the initial condition $\mu_1^*(x_1)=p_1(x_1)$ and the terminal condition
$\phi_{T+1}^*(x_{T+1}, u_{T-n+1}^{T})\triangleq \exp\{-c_{T+1}(x_{T+1})\}$. 
\end{thm}

The optimality condition \eqref{eqoptcond} can be utilized to develop a numerical algorithm to find an optimal solution candidate to  the simplified TERMDP \eqref{eqsimpleprob}.
Unfortunately, it will soon be shown that the optimality condition \eqref{eqoptcond} is only necessary in general. 
Since \eqref{eqoptcond} is a nonlinear condition, it is possible that \eqref{eqoptcond} admits multiple distinct solutions, some of which may correspond to local minima and saddle points of the simplified TERMDP \eqref{eqsimpleprob}.
Theorem~\ref{theooptcond} is closely related to the previously obtained condition in \cite{charalambous2014optimization} and \cite{stavrou2015information}. Theorem~\ref{theooptcond} refines the results of \cite{charalambous2014optimization} and \cite{stavrou2015information} by incorporating the underlying Makovian structure of the simplified TERMDP \eqref{eqsimpleprob}.

\ifdefined\IFAUTART
\begin{figure*}[t]
\rule{\textwidth}{1pt}
\begin{alg} (Forward-Backward Arimoto-Blahut Algorithm)

\label{alg1}
\rule{\textwidth}{1pt}
\begin{itemize}[leftmargin=0ex]
\item[] \textbf{Initialize}
  \begin{flalign*}
  \hspace{5ex}&q_t^{(0)}(u_t|x_{t-m}^t, u_{t-n}^{t-1}) \text{ for } t=1,2, ... , T; \\
  &\phi_{T+1}^{(k)}(x_{T+1-m}^{T+1}, u_{T+1-n}^{T})\triangleq \exp\{-C_{T+1}(x_{T+1})\} \text{ for } k=1,2, ..., K; &&
  \end{flalign*}
\item[] \textbf{For} $k=1, 2, ... , K$ (until convergence) \textbf{do}
\begin{siderules}
	\begin{itemize}[leftmargin=2ex]
	\item[] \textbf{For} $t=1, 2, ... , T-1, T$  (forward path) \textbf{do}
	\begin{siderules}
   \begin{flalign*}
   \hspace{3ex}&\mu_{t+1}^{(k)}(x_{t-m+1}^{t+1}, u_{t-n+1}^t) =\sum\nolimits_{\mathcal{X}_{t-m}}\sum\nolimits_{\mathcal{U}_{t-n}}p_{t+1}(x_{t+1}|x_t,u_t) q_t^{(k-1)}(u_t|x_{t-m}^t, u_{t-n}^{t-1})\mu_t^{(k)}(x_{t-m}^t, u_{t-n}^{t-1}); \\
   &\nu_t^{(k)}(u_t|u_{t-n}^{t-1})=\sum\nolimits_{\mathcal{X}_{t-m}^t} \!\!\! q_t^{(k-1)}(u_t|x_{t-m}^t, u_{t-n}^{t-1})\mu_t^{(k)}(x_{t-m}^t|u_{t-n}^{t-1}); &&
   \end{flalign*}
   \end{siderules}
	\item[] \textbf{End}
	\end{itemize}	
	\begin{itemize}[leftmargin=2ex]
	\item[] \textbf{For} $t=T, T-1, ... , 2, 1$ (backward path) \textbf{do}
\begin{siderules}
\begin{flalign*}
\hspace{3ex}&\rho^{(k)}_t(x_{t-m}^t, u_{t-n}^t)= c_t(x_t,u_t) - \sum_{\calx_{t+1}} p_{t+1}(x_{t+1}|x_t,u_t) \log \phi^{(k)}_{t+1}(x_{t-m+1}^{t+1}, u_{t-n+1}^{t}); \\
&\phi^{(k)}_t(x_{t-m}^t, u_{t-n}^{t-1}) = \sum_{\mathcal{U}_t}\nu^{(k)}_t(u_t|u^{t-1}_{t-n})\exp \left\{-\rho^{(k)}_t(x_{t-m}^t, u_{t-n}^t)\right\};  \\
&q_t^{(k)}(u_t|x_{t-m}^t, u_{t-n}^{t-1}) =\frac{\nu^{(k)}_t(u_t|u^{t-1}_{t-n})\exp \left\{-\rho^{(k)}_t(x_{t-m}^t, u_{t-n}^t)\right\}}{\phi^{(k)}_t(x_{t-m}^t, u_{t-n}^{t-1})};  &&
\end{flalign*}
\end{siderules}
	\item[] \textbf{End}
	\end{itemize}
\end{siderules}
\item[] \textbf{End}
\item[] \textbf{Return} $q_t^{(K)}(u_t|x_{t-m}^t, u_{t-n}^{t-1}), t=1, 2, ... , T$;
\end{itemize}
\end{alg}
\rule{\textwidth}{1pt}
\end{figure*}
\fi

\ifdefined\IFIEEE
\begin{algorithm*}
  \label{alg1}
  \caption{Forward-Backward Arimoto-Blahut Algorithm}
  Initialize:
   \begin{align} 
   &q_t^{(0)}(u_t|x_t, u_{t-n}^{t-1})>0 \text{ for  } t=1,2, ... , T;\label{eqalg1init} \\
  &\phi_{T+1}^{(k)}(x_{T+1}, u_{T+1-n}^{T})\triangleq \exp\{-c_{T+1}(x_{T+1})\} \text{ for } k=1,2, ..., K; \nonumber
  \end{align}
  \For{$k=1, 2, ... , K$ (until convergence)}
  {
   $//$ (Forward path)\;
   \For{$t=1, 2, ... , T$}
   {
   \begin{subequations}
   \begin{align}
   &\mu_{t+1}^{(k)}(x_{t+1}, u_{t-n+1}^t) =\sum\nolimits_{x_t\in\mathcal{X}_t}\sum\nolimits_{u_{t-n}\in\mathcal{U}_{t-n}}p_{t+1}(x_{t+1}|x_t,u_t) q_t^{(k-1)}(u_t|x_t, u_{t-n}^{t-1})\mu_t^{(k)}(x_t, u_{t-n}^{t-1}); \label{eqalg1_1} \\
   &\nu_t^{(k)}(u_t|u_{t-n}^{t-1})=\sum\nolimits_{x_t\in\mathcal{X}_t}  q_t^{(k-1)}(u_t|x_t, u_{t-n}^{t-1})\mu_t^{(k)}(x_t|u_{t-n}^{t-1}); \label{eqalg1_2}
   \end{align}
   \end{subequations}

   }
   $//$ (Backward path)\;
   \For{$t=T, T-1, ... , 1$}
   {
   \begin{subequations}
   \label{eqabback}
   \begin{align}
   &\rho^{(k)}_t(x_t, u_{t-n}^t)= c_t(x_t,u_t) - \sum\nolimits_{x_{t+1}\in\calx_{t+1}} p_{t+1}(x_{t+1}|x_t,u_t) \log \phi^{(k)}_{t+1}(x_{t+1}, u_{t-n+1}^{t})  ;  \label{eqalg1_3}\\
   &\phi^{(k)}_t(x_t, u_{t-n}^{t-1}) =  \sum\nolimits_{u_t \in \mathcal{U}_t}\nu^{(k)}_t(u_t|u^{t-1}_{t-n})\exp \left\{-\rho^{(k)}_t(x_t, u_{t-n}^t)\right\}; \label{eqalg1_4}\\
   &q_t^{(k)}(u_t|x_t, u_{t-n}^{t-1}) =\tfrac{\nu^{(k)}_t(u_t|u^{t-1}_{t-n})\exp \left\{-\rho^{(k)}_t(x_t, u_{t-n}^t)\right\}}{\phi^{(k)}_t(x_t, u_{t-n}^{t-1})};\label{eqalg1_5}
   \end{align}
   \end{subequations}
   }
  }
  Return $q_t^{(K)}(u_t|x_t, u_{t-n}^{t-1})$\;
\end{algorithm*}
\fi

\subsection{Forward-Backward Arimoto-Blahut Algorithm}
As the second contribution, we propose an iterative algorithm to solve \eqref{eqoptcond} numerically. To this end, we classify the five equations into two groups. Equations \eqref{eqoptcond1} and \eqref{eqoptcond2} form the first group (characterizing  variables $\mu^*$ and $\nu^*$), and equations \eqref{eqoptcond3}-\eqref{eqoptcond5} form the second group (characterizing variables $\rho^*$, $\phi^*$, and  $q^*$).
Observe that if the variables $(\rho^*, \phi^*, q^*)$ are known, then the first set of equations, which can be viewed as the Kolmogorov forward equation, can be solved forward in time to compute $(\mu^*, \nu^*)$.
Conversely, if the variables $(\mu^*, \nu^*)$ are known, then the second set of equations, which can be viewed as the Bellman backward equation, can be solved backward in time to compute  $(\rho^*, \phi^*, q^*)$.
Hence, to compute these unknowns simultaneously, the following boot-strapping method is natural: first, the forward computation is performed using the current best guess of the second set of unknowns, and then the backward computation is performed using the updated guess of the first set of unknowns. The forward-backward iteration is repeated {\upd sufficiently many times}. The proposed algorithm is summarized in Algorithm~\ref{alg1}.
Notice that Algorithm~\ref{alg1} is a generalization of the standard Arimoto-Blahut algorithm \eqref{eqarimotoblahut}, which can be recovered as a special case with $T=1$, $m=n=0$ and $C_{T+1}(\cdot)=0$.

{\updnew
Clearly, solutions of \eqref{eqoptcond} are fixed points of Algorithm~\ref{alg1}. 
The second main result of this paper is stated as follows.}
{\upd
\begin{thm}
\label{theoconvergefbaba}
For each initial condition \eqref{eqalg1init},  the sequence $q^{(k)}$ generated by Algorithm~\ref{alg1} has a limit point. Moreover, every limit point $q^{*}$  satisfies \eqref{eqoptcond}.
\end{thm}

Theorems \ref{theooptcond} and \ref{theoconvergefbaba} together imply that every limit point of Algorithm~\ref{alg1} is an optimal solution candidate for TERMDP.
The following remarks clarify some limitations of our main results, which require further analysis in future work.
\begin{remark}
\label{rem1}
\begin{enumerate}[leftmargin=*]
\item  Since \eqref{eqoptcond} is only a necessary optimality condition, a limit point may not be a local minimum of the given TERMDP. Section~\ref{secnonconvex} below presents an example in which one of the limit points is a saddle point.
\item Theorem~\ref{theoconvergefbaba} does not guarantee the existence of $\lim_{k\rightarrow \infty} q^{(k)}$ (e.g., the sequence may oscillate between distinct points). However, we were not able to construct such a counterexample.
\item In the classical Arimoto-Blahut algorithm, there is a known stopping criterion by which suboptimality of the iteration is estimated effectively (e.g., \cite[Fig. 3]{blahut1972computation}). Currently, such a stopping criterion is not known for Algorithm~\ref{alg1}.
\end{enumerate}
\end{remark}

The number of arithmetic operations to perform a single forward-backward path in Algorithm~\ref{alg1}
is estimated as $O(T|\calx|^2|\calu|^{n+1})$. Notice that it is linear in $T$, grows exponentially with $n$, and does not depend on $m$.

}

\subsection{Nonconvexity}
\label{secnonconvex}
Due to the nonconvexity of the value functions in \eqref{eqhjb2}, the optimality condition \eqref{eqoptcond} is only necessary in general.
In fact, depending on the initial condition, Algorithm~\ref{alg1} can converge to different stationary points corresponding to local minima and saddle points of the considered TERMDP \eqref{eqsimpleprob}.
To demonstrate this, consider a simple problem instance of the TERMDP \eqref{eqsimpleprob} with $T=2$, $n=0$, $\mathcal{X}=\{0, 1\}$, $\mathcal{U}=\{0, 1\}$,
\[
c_t(x_t, u_t)=\begin{cases}
0 & \text{ if } u_t=x_t \\
1 & \text{ if } u_t\neq x_t
\end{cases}
\]
for $t=1, 2$, and $c_3(x_3)\equiv 0$.
The initial state distribution is assumed to be 
$\mu_1(x_1=0)=\mu_1(x_1=1)=0.5$, and the state transitions are deterministic in that
\[
p_{t+1}(x_{t+1}| x_t, u_t)=\begin{cases}
1 & \text{ if } x_{t+1}=u_t \\
0 & \text{ if } x_{t+1}\neq u_t .
\end{cases}
\]
To see that this problem has multiple distinct local minima, we solve the Bellman equation \eqref{eqhjb2} numerically by griding the space of probability distributions.
Specifically, at $t=2$, the value function $V_2(\mu_2(x_2))$ is computed by solving the minimization 
\[
V_2(\mu_2(x_2))=\min_{q_2(u_2|x_2)} \{\mathbb{E}c_2({\upd X_2, U_2})+I(X_2; U_2)\}.
\]
Since $\mu_2(x_2)$ is an element of a single-dimensional probability simplex, it can be parameterized as $\mu_2(x_2=0)=\lambda$ and  $\mu_2(x_2=1)=1-\lambda$ with $\lambda\in[0,1]$.
For each fixed $\lambda$, the minimization above is solved by the standard Arimoto-Blahut iteration:
\begin{align*}
\nu_2^{(k)}(u_2)&=\sum_{x_2\in\{0, 1\}}\mu_2(x_2)q_2^{(k-1)}(u_2|x_2) \\
\phi_2^{(k)}(x_2)&=\sum_{u_2\in\{0, 1\}}\nu_2^{(k)}(u_2)\exp\{-c_2(x_2, u_2)\}  \\
q_2^{(k)}(u_2|x_2)&=\frac{\nu_2^{(k)}(u_2)\exp\{-c_2(x_2, u_2)\}}{\phi_2^{(k)}(x_2)}.
\end{align*}
After the convergence, the value function is computed as
\[
V_2(\mu_2(x_2))=-\sum_{x_2\in\{0, 1\}} \mu_2(x_2)\log \phi_2^{(k)}(x_2).
\]
Fig.~\ref{fig:contour} (Left) shows  $V_2(\mu_2(x_2))$ as a function of $\lambda$. It is clearly nonconvex. After  $V_2(\mu_2(x_2))$ is obtained, the Bellman equation at time $t=1$ can be evaluated as
\begin{equation}
\label{eqvalf1}
V_1(\mu_1(x_1))=\min_{q_1(u_1|x_1)} \{\mathbb{E}c_1({\upd X_1, U_1})+I(X_1; U_1)+V_2(\mu_2(x_2))\}.
\end{equation}
Due to the nonconvexity of $V_2(\mu_2(x_2))$, the objective function in the minimization \eqref{eqvalf1} is a nonconvex function of $q_1(u_1|x_1)$. Fig.~\ref{fig:contour} (Right) shows the objective function in \eqref{eqvalf1} plotted as a function of $q_1$ parameterized by $\theta_0$ and $\theta_1$:
\begin{align*}
q_1(u_1=0|x_1=0)&=\theta_0 \\
q_1(u_1=1|x_1=0)&=1-\theta_0 \\
q_1(u_1=0|x_1=1)&=\theta_1 \\
q_1(u_1=1|x_1=1)&=1-\theta_1. 
\end{align*}
Clearly, it is a nonconvex function, admitting two local minima (A and C) and a saddle point (B). Each of them is a fixed point of Algorithm~\ref{alg1}.

\ifdefined\IFONECOL

\else
\begin{figure}[t]
    \centering
    \includegraphics[width=\columnwidth]{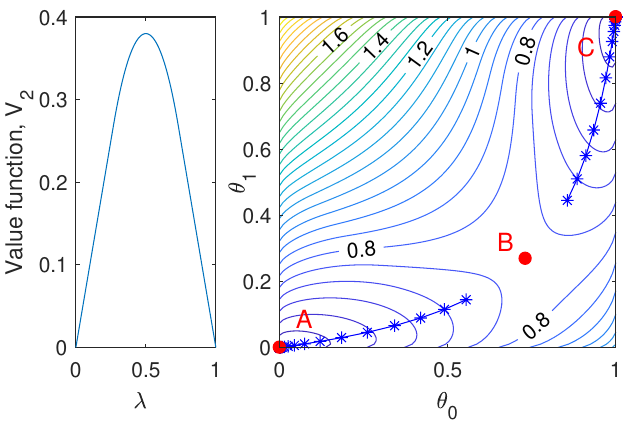}
    \caption{Left: The value function $V_2$. Right: Contour plot of the objective function $\mathbb{E}c_1(x_1, u_1)+I(X_1; U_1)+V_2(\mu_2(x_2))$ in \eqref{eqvalf1} as a function of $\theta_0$ and $\theta_1$. Stationary points A and B are local minima, whereas C is a saddle point. Two sample trajectories of the proposed forward-backward Arimoto-Blahut algorithm (Algorithm~\ref{alg1}) started with different initial conditions are also shown. It can be shown that A, B and C are all fixed points of Algorithm~\ref{alg1}.}
    \label{fig:contour}
\end{figure}
\fi

\section{Derivation of Main Results}
\label{secderivation}
{\updnew This section provides technical details of the proofs of Theorems~\ref{theooptcond} and \ref{theoconvergefbaba}.
}
\subsection{Preparation}
The first step is to rewrite the objective function in \eqref{eqsimpleprob} as an explicit function of $q^T$.
For each $t=1, 2, ... , T$, let $\mu_t(x_t|u_{t-n}^{t-1})$ be the conditional distribution obtained from $\mu_t(x_t, u_{t-n}^{t-1})$ whenever $\mu_t(u_{t-n}^{t-1})>0$. 
Define the conditional distribution $\nu_t$ by
\begin{equation}
\nu_t(u_t|u_{t-n}^{t-1})=\sum\nolimits_{x_t\in\mathcal{X}_t}q_t(u_t|x_t, u_{t-n}^{t-1})\mu_t(x_t|u_{t-n}^{t-1}).\label{eqdefnu}
\end{equation}
When $\mu_t(u_{t-n}^{t-1})=0$, $\nu_t$ is defined to be the uniform distribution on $\mathcal{U}_t$.
For each $t=1, 2, ... , T$, we consider $\nu_t$ and $q_t$ as elements of Euclidean spaces, i.e.,
\begin{subequations}
\label{eqnuqeuc}
\begin{align}
\nu_t(u_t|u_{t-n}^{t-1})&\in \mathbb{R}^{|\calu_{t-n}|\times \cdots \times |\calu_t|} \label{eqnuqeuc1} \\
q_t(u_t|x_t, u_{t-n}^{t-1})&\in \mathbb{R}^{|\calu_{t-n}|\times \cdots \times |\calu_t|\times |\calx_t|}. \label{eqnuqeuc2}
\end{align}
\end{subequations}
Since $\nu_t$ and $q_t$ are conditional probability distributions, they are entry-wise non-negative  (denoted by $\nu_t\geq 0$ and $q_t\geq 0$) and
\begin{subequations}
\label{eqnuqconst}
\begin{align}
&\sum_{u_t\in \calu_t} \nu_t(u_t|u_{t-n}^{t-1})=1 \; \forall u_{t-n}^{t-1}\in\calu_{t-n}^{t-1}  \label{eqnuqconst1}\\
&\sum_{u_t\in \calu_t} q_t(u_t|x_{t}, u_{t-n}^{t-1})=1 \; \forall (x_{t}, u_{t-n}^{t-1}) \in \calx_{t} \times \calu_{t-n}^{t-1}. \label{eqnuqconst2}
\end{align}
\end{subequations}
Thus, the feasibility sets for $\nu^T$ and $q^T$ are 
\begin{align*}
X_\nu=\{\nu^T: \text{\eqref{eqnuqeuc1}, \eqref{eqnuqconst1} and } \nu_t\geq 0 \text{ for every } t=1, 2, ... , T\}; \\
X_q=\{q^T: \text{\eqref{eqnuqeuc2}, \eqref{eqnuqconst2} and } q_t\geq 0 \text{ for every } t=1, 2, ... , T\}.
\end{align*}
Using $\mu_t$, $\nu_t$ and $q_t$, the stage-wise cost in \eqref{eqsimpleprob} can be written as 
\begin{align*}
\ell_t(\mu_t, \nu_t, q_t) &\triangleq \mathbb{E}c_t({\upd X_t, U_t})+I(X_t;U_t|U_{t-n}^{t-1})\\
&=\sum_{x_k\in\calx_k}\sum_{u_{k-n}^k\in\calu_{k-n}^k} \mu_k(x_k, u_{k-n}^{k-1}) q_k(u_k|x_k, u_{k-n}^{k-1}) \\
&\hspace{8ex}\times \left( \!\log \!\frac{q_k(u_k|x_k,u_{k-n}^{k-1})}{\nu_k(u_k|u_{k-n}^{k-1})} \!+\! c_k(x_k,u_k) \!\right)
\end{align*}
for $t=1, 2, ... , T$ and 
\begin{align*}
\ell_{T+1}(\mu_{T+1})&\triangleq \mathbb{E} c_{T+1}({\upd X_{T+1}}) \\
&=\sum\nolimits_{x_{T+1}\in\mathcal{X}_{T+1}}\mu_{T+1}(x_{T+1})c_{T+1}(x_{T+1})
\end{align*}
for the final time step.
Thus, the objective function in \eqref{eqsimpleprob} is 
\begin{equation}
\label{eqfvuq}
f(\nu^T, q^T)=\sum\nolimits_{t=1}^T \ell_t(\mu_t, \nu_t, q_t)+\ell_{T+1}(\mu_{T+1}).
\end{equation}
{\updnew 
Notice that we consider \eqref{eqfvuq} as a function of $\nu^T$ and $q^T$, which turns out to be convenient in order to view Algorithm~\ref{alg1} as a two-block BCD algorithm. Our problem is to minimize $f(\nu^T, q^T)$ over $\mathcal{X}_\nu\times\mathcal{X}_q$ subject to the equality  constraint \eqref{eqdefnu}.
However, the following result states that \eqref{eqdefnu} will be automatically satisfied  by a coordinate-wise minimizer $\nu^T$ of $f(\nu^T, q^T)$ for a fixed $q^T$. 
}
\begin{lem}\cite[Theorem 4(b)]{blahut1972computation}
\label{lemblahut}
Let $q^T\in X_q$ be fixed. Then
\begin{subequations}
\begin{align}
\min & \hspace{2ex} f(\nu^T, q^T) \\
\text{s.t. } & \hspace{2ex} \nu^T \in X_\nu
\end{align}
\end{subequations}
is a convex optimization problem with respect to $\nu^T$, and an optimal solution is given by \eqref{eqdefnu}.
\end{lem}

{\updnew
Therefore, the simplified TERMDP \eqref{eqsimpleprob} can be written as
\begin{subequations}
\label{eqsimpleprob2}
\begin{align}
\min & \hspace{2ex} f(\nu^T, q^T) \\
\text{s.t. } & \hspace{2ex} \nu^T\in X_\nu, q^T\in X_q.
\end{align}
\end{subequations}
(Including the equality constraint \eqref{eqdefnu} in \eqref{eqsimpleprob2} does not alter the the optimal solution because of Lemma~\ref{lemblahut}.)
Notice that if $q^{T*}$ is a local minimizer for \eqref{eqsimpleprob}, then there exists a vector $\nu^{T*}\in\mathcal{X}_\nu$ such that $(\nu^{T*}, q^{T*})$ is a local minimizer for \eqref{eqsimpleprob2}. Moreover, a local minimizer $(\nu^{T*}, q^{T*})$ is necessarily a coordinate-wise local minimizer. By coordinate-wise convexity (which will be shown below) of $f(\nu^T, q^T)$, we have:
\begin{subequations}
\label{eqblockoptimal}
\begin{align}
\nu^{T*}&\in \argmin\nolimits_{\nu^T \in X_\nu} f(\nu^T, q^{T*})  \label{eqblockoptimal1} \\
q^{T*}&\in \argmin\nolimits_{q^T \in X_q} f(\nu^{T*}, q^T) \label{eqblockoptimal2}.
\end{align}
\end{subequations}
Therefore, \eqref{eqblockoptimal} is a necessary condition for $q^{T*}$ to be a local minimizer for the simplified TERMDP \eqref{eqsimpleprob}.
In the next subsection, we further show that \eqref{eqblockoptimal} implies the condition \eqref{eqoptcond}, which thus shows that  \eqref{eqoptcond} is a necessary optimality condition for the simplified TERMDP \eqref{eqsimpleprob}. This argument outlines the proof of Theorem~\ref{theooptcond}, which will be detailed in Section~\ref{secproofmainresults}.
}


{\upd
\subsection{Analysis of Algorithm~\ref{alg1}}
For the analysis of Algorithm~\ref{alg1}, a key observation is that it is the two-block BCD  algorithm  applied to \eqref{eqblockoptimal}:
\begin{subequations}
\label{eqbcdfb}
\begin{align}
\nu^{T(k)}&=\argmin\nolimits_{\nu^T \in X_\nu} f(\nu^T,q^{T(k-1)})  \label{eqbcdfb1} \\
q^{T(k)}&=\argmin\nolimits_{q^T \in X_q} f(\nu^{T(k)}, q^T).  \label{eqbcdfb2}
\end{align}
\end{subequations}
To prove that the backward path \eqref{eqabback} is equivalent to the coordinate-wise minimization  \eqref{eqbcdfb2}, we remark that} the minimization problem \eqref{eqbcdfb2} can be viewed as an optimal control problem with respect to the control actions $q^T=(q_1, ... , q_T)$. The next lemma essentially shows that the backward path \eqref{eqabback}  is solving \eqref{eqbcdfb2}  by  backward dynamic programming.
\begin{lem}
\label{lemblockq}
Suppose
\begin{align*}
&(\nu_1^{(k)}, ... , \nu_T^{(k)}) \in X_\nu, \\
&(q_1^{(k-1)}, ... , q_{\tau-1}^{(k-1)}, q_\tau, q_{\tau+1}^{(k)}, ... , q_T^{(k)}) \in X_q,
\end{align*}
and $\{\mu_{t+1}\}_{t=1}^T$ is the sequence of probability measures generated by
$(q_1^{(k-1)}, ... , q_{\tau-1}^{(k-1)}, q_\tau, q_{\tau+1}^{(k)}, ... , q_T^{(k)})$ via  \eqref{eqmustate}.
Let $\rho_\tau^{(k)}, \phi_\tau^{(k)}$ and $q_\tau^{(k)}$ be the parameters obtained by computing \eqref{eqabback} backward in time for $t=T, ... , \tau$. Then, for each $\tau=T, T-1, ... ,1$, the following statements hold:
\begin{itemize}
\item[(a)] The function
\[
f(\nu_1^{(k)}, ... , \nu_T^{(k)}, q_1^{(k-1)}, ... , q_{\tau-1}^{(k-1)}, q_\tau, q_{\tau+1}^{(k)}, ... , q_T^{(k)})
\]
is convex in $q_\tau \geq 0$, and any global minimizer $q_\tau^\circ$ satisfies $q_\tau^\circ=q_\tau^{(k)}$ almost everywhere with respect to $\mu_\tau(x_\tau, u_{\tau-n}^{\tau-1})$.
\item[(b)] The cost-to-go function under the policy $\{q_t^{(k)}\}_{t=\tau}^T$ is linear in $\mu_\tau$:
\begin{align*}
&\sum\nolimits_{t=\tau}^T \ell_t(\mu_t, \nu_t^{(k)}, q_t^{(k)})+\ell_{T+1}(\mu_{T+1})\\
&=-\sum_{x_\tau\in\mathcal{X}_\tau}\sum_{u_{\tau-n}^{\tau-1}\in\mathcal{U}_{\tau-n}^{\tau-1}} \mu_\tau(x_\tau, u_{\tau-n}^{\tau-1})\log \phi_\tau^{(k)}(x_\tau, u_{\tau-n}^{\tau-1}).
\end{align*}
\end{itemize}
\end{lem}
\begin{pf}
The proof is by backward induction. For the time step $T$, we have
\begin{align}
&f(\nu_1^{(k)}, ... , \nu_T^{(k)}, q_1^{(k-1)}, ... , q_{T-1}^{(k-1)}, q_T) \nonumber \\
&= \sum_{x_T\in\calx_T}\sum_{u_{T-n}^T\in\calu_{T-n}^T} \mu_T(x_T, u_{T-n}^{T-1})  q_T(u_T|x_T, u_{T-n}^{T-1}) \nonumber\\
& \times \!\left(\log\frac{q_T(u_T|x_T, u_{T-n}^{T-1})}{\nu^{(k)}_T(u_T| u_{T-n}^{T-1})}\!+\!\rho_T^{(k)}(x_T, u_{T-n}^T) \right)+\text{const.}  \label{eqqtlast}
\end{align}
{\updnew where $\mu_T$, $\nu_T^{(k)}, \rho_T^{(k)}$ and ``const.'' do not depend on $q_T$. Convexity of \eqref{eqqtlast} in $q_T$ is clear, as $q_T\log q_T$ is a convex function in $q_T$. Moreover, minimization of \eqref{eqqtlast} in terms of $q_T$ has the same structure as the minimization problem \eqref{eqrdprob} in terms of $q(u|x)$. Therefore, it follows from Proposition~\ref{propext} that the minimizer $q_T^\circ$ for \eqref{eqqtlast} is given by $q_T^\circ=q_T^{(k)}$.}
This establishes (a) for the time step $\tau=T$. The statement (b) for $\tau=T$ can be directly shown by substituting the expression of $q_t^{(k)}$ given by
\eqref{eqalg1_5} with {\updnew  $t=T$}
into \eqref{eqqtlast}:
\begin{subequations}
\label{eqind_b}
\begin{align}
&\ell_T(\mu_T, \nu_T^{(k)}, q_T^{(k)})+\ell_{T+1}(\mu_{T+1})\nonumber \\
&= \sum_{x_T\in\calx_T}\sum_{u_{T-n}^T\in\calu_{T-n}^T} \mu_T(x_T, u_{T-n}^{T-1})  q_T^{(k)}(u_T|x_T, u_{T-n}^{T-1}) \nonumber\\
& \hspace{3ex}\times \left(\log\frac{q_T^{(k)}(u_T|x_T, u_{T-n}^{T-1})}{\nu^{(k)}_T(u_T| u_{T-n}^{T-1})}+\rho_T^{(k)}(x_T, u_{T-n}^T) \right)  \\
&= \sum_{x_T\in\calx_T}\sum_{u_{T-n}^T\in\calu_{T-n}^T} \mu_T(x_T, u_{T-n}^{T-1})  q_T^{(k)}(u_T|x_T, u_{T-n}^{T-1}) \nonumber\\
& \hspace{3ex}\times \left(-\log \phi_T^{(k)}(x_T, u_{T-n}^{T-1})\right) \\
&= -\sum_{x_T\in\calx_T}\sum_{u_{T-n}^{T-1}\in\calu_{T-n}^{T-1}} \mu_T(x_T, u_{T-n}^{T-1})  \log \phi_T^{(k)}(x_T, u_{T-n}^{T-1}) \nonumber \\
& \hspace{3ex}\times \underbrace{\sum_{u_T\in\calu_T}q_T^{(k)}(u_T|x_T, u_{T-n}^{T-1})}_{=1}.
\end{align}
\end{subequations}

To complete the proof, {\updnew we show that if (b) holds for the time step $\tau+1$, then both (a) and (b) hold for the time step $\tau$.} Since (b) is hypothesized for $\tau+1$, using $\rho_\tau^{(k)}$, it is possible to write
\begin{align}
&f(\nu_t^{(k)}, ... , \nu_T^{(k)}, q_1^{(k-1)}, ... , q_{\tau-1}^{(k-1)}, q_\tau, q_{\tau+1}^{(k)}, ... , q_T^{(k)}) \nonumber \\
&= \sum_{x_\tau\in\calx_\tau}\sum_{u_{\tau-n}^\tau\in\calu_{\tau-n}^\tau} \mu_\tau(x_\tau, u_{\tau-n}^{\tau-1})  q_\tau(u_\tau|x_\tau, u_{\tau-n}^{\tau-1}) \nonumber\\
& \times \left(\log\frac{q_\tau(u_\tau|x_\tau, u_{\tau-n}^{\tau-1})}{\nu^{(k)}_\tau(u_\tau| u_{\tau-n}^{\tau-1})}+\rho_\tau^{(k)}(x_\tau, u_{\tau-n}^\tau) \right)+\text{const.}  \label{eqqtmid2}
\end{align}
{\updnew where $\mu_\tau$, $\nu_\tau^{(k)}$, $\rho_\tau^{(k)}$ and ``const.'' do not depend on $q_\tau$. Namely, \eqref{eqqtmid2} is convex in $q_\tau$.}
Proposition~\ref{propext} is applicable once again to conclude that a minimizer coincides with $q_\tau^{(k)}$ almost everywhere with respect to $\mu_\tau(x_\tau, u_{\tau-n}^{\tau-1})$.
Hence, (a) is established for the time step $\tau$. The statement (b) for $\tau$ can be shown by a direct substitution. The details are similar to \eqref{eqind_b}.
\end{pf}

\subsection{Proof of main results}
\label{secproofmainresults}
Based on the observations so far, the main theorems in this paper are established as follows.

\textbf{Proof of Theorem~\ref{theooptcond}:}
For a given $\{q_t^*\}_{t=1}^T$, we define $\{\mu_t^*\}_{t=1}^T$ via \eqref{eqoptcond1} and hence \eqref{eqoptcond1} is automatically satisfied. By Lemma~\ref{lemblahut}, for any locally optimal solution $\{q_t^*\}_{t=1}^T$ to \eqref{eqsimpleprob}, there exist variables $\{\nu_t^*\}_{t=1}^T$ satisfying \eqref{eqoptcond2}, such that $\{\nu_t^*, q_t^*\}_{t=1}^T$ is a locally optimal solution to \eqref{eqsimpleprob2}.
Since local optimality implies coordinate-wise local optimality, for each $t=1, 2, ... , T$, $q_t^*$ is a local minimizer of
\begin{equation}
\label{eqfnuqstar}
f(\nu_1^*, ... , \nu_T^*, q_T^*, ... , q_{t+1}^*, q_t, q_{t-1}^*, ... , q_1^*).
\end{equation}
However, Lemma~\ref{lemblockq} is applicable (with $\nu_\tau^{(k)}=\nu_\tau^*$ for $\tau=1, ... , T$, $q_\tau^{(k)}=q_\tau^*$ for $\tau=t+1, ... , T$ and $q_\tau^{(k-1)}=q_\tau^*$ for $\tau=1, ... , t-1$) to conclude that \eqref{eqfnuqstar} is convex in $q_t\geq 0$ and hence
\[
q_t^* \in \argmin_{q_t\geq 0} f(\nu_1^*, ... , \nu_T^*, q_T^*, ... , q_{t+1}^*, q_t, q_{t-1}^*, ... , q_1^*).
\]
Moreover, Lemma~\ref{lemblockq} (a) also implies that if the parameters $\{\rho_\tau^*, \phi_\tau^*\}_{\tau=t}^T$ are calculated by \eqref{eqoptcond3}-\eqref{eqoptcond4} backward in time, then any global minimizer
\[
q_t^\circ \in \argmin_{q_t\geq 0} f(\nu_1^*, ... , \nu_T^*, q_T^*, ... , q_{t+1}^*, q_t, q_{t-1}^*, ... , q_1^*)
\]
satisfies
\[
q_t^\circ=\frac{\nu^*_t(u_t|u^{t-1}_{t-n})\exp \left\{-\rho^*_t(x_t, u_{t-n}^t)\right\}}{\phi^*_t(x_t, u_{t-n}^{t-1})}
\]
$\mu_t$-almost everywhere. Hence \eqref{eqoptcond5} must hold. \hfill \QED

\textbf{Proof of Theorem~\ref{theoconvergefbaba}:}
{\upd
The fact that the sequence $q^{(k)}$ generated by Algorithm~\ref{alg1} has a limit point follows from  the fact that $X_q$ is a compact set.

To show that every limit point of the sequence $(\nu^{(k)}, q^{(k)})$ generated by Algorithm~\ref{alg1} is a stationary point for \eqref{eqsimpleprob2}, observe that
\begin{itemize}
\item[(a)] The update rule \eqref{eqalg1_2} is equivalent to \eqref{eqbcdfb1}; and
\item[(b)] The update rule  \eqref{eqalg1_5} is equivalent to \eqref{eqbcdfb2}.
\end{itemize}
The fact (a) follows from Lemma~\ref{lemblahut} and (b) follows from Lemma~\ref{lemblockq}.
Therefore, Algorithm~\ref{alg1} is equivalent to the two-block BCD algorithm to which Lemma~\ref{lembcd} is applicable.

Finally, we claim that every stationary point for \eqref{eqsimpleprob2} satisfies \eqref{eqoptcond}.
To see this, notice that every stationary point $(\nu^{T*}, q^{T*})$ is a coordinate-wise stationary point. Since $f(\nu^{T}, q^{T})$ is coordinate-wise convex, $(\nu^{T*}, q^{T*})$ is a coordinate-wise minimizer, i.e., \eqref{eqblockoptimal} holds. By Lemma~\ref{lemblahut}, conditions \eqref{eqoptcond1} and \eqref{eqoptcond2} are necessary for \eqref{eqblockoptimal1}.
By Lemma~\ref{lemblockq}  (with $\nu_\tau^{(k)}=\nu_\tau^*$ for $\tau=1, ... , T$, $q_\tau^{(k)}=q_\tau^*$ for $\tau=t+1, ... , T$ and $q_\tau^{(k-1)}=q_\tau^*$ for $\tau=1, ... , t-1$),  conditions  \eqref{eqoptcond3}-\eqref{eqoptcond5} are necessary for \eqref{eqblockoptimal2}.
\hfill \QED
}

\section{Interpretations}
\label{secinterpret}

In this section, we discuss two applications of TERMDP in engineering and scientific contexts in which  transfer entropy plays central roles.

\subsection{Networked Control Systems}
\label{secncs}

The first application is the analysis of networked control systems, where the sensor data is transmitted to the controller over a rate-limited communication channel.
Fig.~\ref{fig:channel} shows a discrete-time, finite-horizon MDP setup in which a decision policy must be realized by a joint design of encoder and decoder, together with an appropriate codebook for discrete noiseless channel.
Most generally, assume that an encoder is a stochastic kernel $e_t(w_t|x^t, w^{t-1})$ and a decoder is a stochastic kernel $d_t(u_t|w^t,u^{t-1})$.
At each time step, a codeword $w_t$ is chosen from a codebook $\mathcal{W}_t$ such that $|\mathcal{W}_t|=2^{R_t}$.
We refer to $R=\sum_{t=1}^T R_t$ as the \emph{rate} of communication.
The next proposition claims that the rate of communication in Fig.~\ref{fig:channel} is fundamentally lower bounded by the directed information.
\begin{prop}
\label{proplbdi}
Let the encoder and the decoder be any stochastic kernels of the form $e_t(w_t|x^t,w^{t-1})$ and $d_t(u_t|w^t, u^{t-1})$. Then $R\log 2\geq I(X^T\rightarrow U^T)$.
\end{prop}
\begin{pf}
Note that
\begin{align*}
R\log 2\! &=\!\sum\nolimits_{t=1}^T R_t \log 2 \\
&\geq \!\sum\nolimits_{t=1}^T \!H(W_t) \\
&\geq \!\sum\nolimits_{t=1}^T \!H(W_t|W^{t-1},U^{t-1}) \\
&\geq \!\sum\nolimits_{t=1}^T \!H(W_t|W^{t-1}\!\!,U^{t-1})\!-\!H(W_t|X^t\!\!,W^{t-1}\!\!,U^{t-1}) \\
&=\!\sum\nolimits_{t=1}^T \!I(X^t;W_t|W^{t-1},U^{t-1}) \\
&\triangleq I(X^T\rightarrow W^T \| U^{T-1}).
\end{align*}
The first inequality is due to the fact that entropy of a discrete random variable cannot be greater than its log-cardinality. Notice that a factor $\log 2$ appears since we are using the natural logarithm in this paper. The second inequality holds because conditioning reduces entropy. The third inequality follows since entropy is nonnegative. The last quantity is known as the \emph{causally conditioned directed information}  \cite{kramer2003capacity}. The feedback data-processing inequality \cite{tanaka2015lqg}
\[ I(X^T\rightarrow U^T) \leq I(X^T \rightarrow W^T\| U^{T-1})\]
is applicable to complete the proof.
\ifdefined\IFAUTART
\hfill\ensuremath{\square}
\fi
\end{pf}

 Proposition~\ref{proplbdi} provides a fundamental performance limitation of a communication system when both encoder and decoder have full memories of the past. However, it is also meaningful to consider restricted scenarios in which the encoder and decoder have limited memories. For instance:
 \begin{itemize}
 \item[(A)] The encoder stochastic kernel is of the form $e_t(w_t|x^t_{t-m})$ and the decoder stochastic kernel is of the form $d_t(u_t|w_t, u_{t-n}^t)$; or
 \item[(B)] The encoder stochastic kernel is $e_t(w_t|x^t_{t-m}, u_{t-n}^{t-1})$ and the decoder is a deterministic function $u_t=d_t(w_t)$. The encoder has an access to the past control inputs $u_{t-n}^{t-1}$ since they are predictable from the past $w_{t-n}^{t-1}$ because the decoder is a deterministic map.
 \end{itemize}
 The next proposition shows that the transfer entropy of degree $(m,n)$ provides a tighter lower bound in these cases.
\begin{prop}
\label{propcom2}
Suppose that the encoder and the decoder have structures specified by (A) or (B) above. Then
\[ R\log 2 \geq I_{m,n}(X^T\rightarrow U^T). \]
\end{prop}
\begin{pf}
\ifdefined\LONGVERSION
See Appendix~\ref{app1}.
\else
See [xx, Appendix D].
\fi
\end{pf}

By solving the TERMDP \eqref{eqmainprob} with different $\beta\geq 0$, one can draw a trade-off curve between $J(X^{T+1}, U^T)$ and $I(X^T\rightarrow U^T)$. Proposition~\ref{propcom2} means that this trade-off curve shows a fundamental limitation of the achievable control performance under the given data rate.

The tightness of the lower bounds provided by Propositions~\ref{proplbdi} and \ref{propcom2} (i.e., whether it is possible to construct an encoder-decoder pair such that the data rate matches its lower bound while satisfying the desired control performance) is the natural next question. 
In the LQG control setup, this question has been studied in \cite{silva2013characterization, extendedversion,kostina2019rate,1701.06368}. In these references, it is shown that the conservativeness of the lower bound provided by Proposition~\ref{proplbdi} is no greater than a small constant.
Beyond the LQG setting, the question is currently wide open.

\ifdefined\IFONECOL
\begin{figure}[t]
    \centering
    \includegraphics[width=0.6\columnwidth]{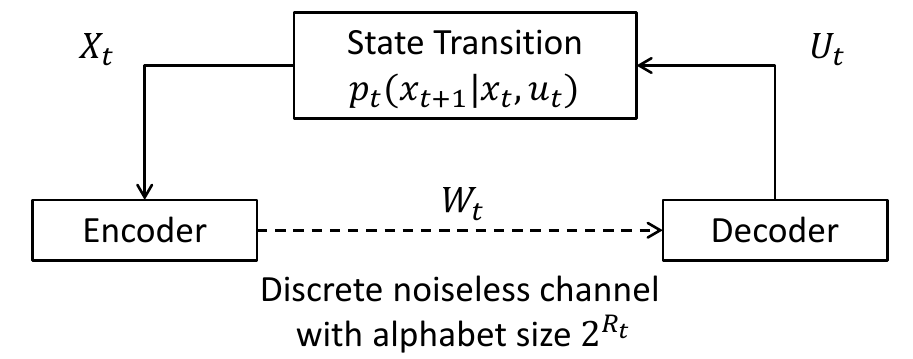}
    \caption{MDP over discrete noiseless channel.}
    \label{fig:channel}
\end{figure}
\else
\begin{figure}[t]
    \centering
    \includegraphics[width=0.8\columnwidth]{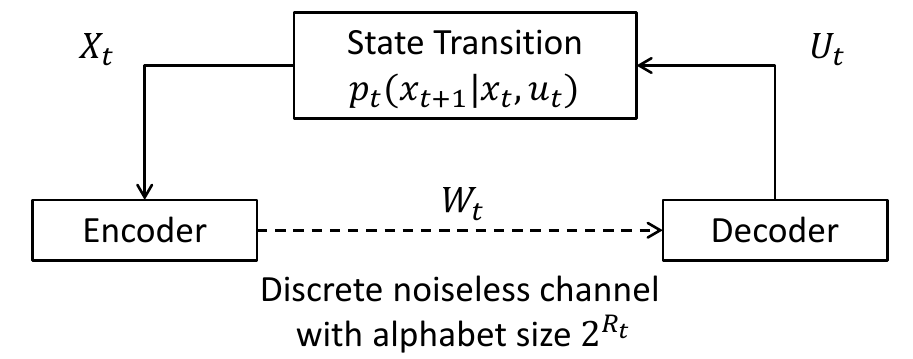}
    \caption{MDP over discrete noiseless channel.}
    \label{fig:channel}
\end{figure}
\fi

\subsection{Maxwell's demon}
\label{secmaxwell}
Maxwell's demon is a physical device that can seemingly violate the second law of thermodynamics, which turns out to be a prototypical thought-experiment that
connects statistical physics and information theory \cite{parrondo2015thermodynamics}.
One of the simplest forms of Maxwell's demon is a device called the Szilard engine. Below, we introduce an application of TERMDP to the analysis of the efficiency of a generalized Szilard engine extracting work at a non-zero rate (in contrast to the common assumption that the engine is operated infinitely slowly).

Consider a single-molecule gas trapped in a box (``engine'') that is immersed in a thermal bath of temperature $T_0$ (Fig.~\ref{fig:szilard}). The state of the engine at time $t$ is represented by the position and the velocity of the molecule, which is denoted by $X_t\in \calx$. Assume that the state space is divided into finite cells so that $\calx$ is a finite set. Also, assume that the evolution of $X_t$ is described by a discrete-time random process.

At each time step $t=0,1, ... , T-1$, suppose that one of the following three possible control actions $U_t$ can be applied: (i) insert a weight-less barrier into the middle of the engine box and move it to the left at a constant velocity $v$ for a unit time, (ii) insert a barrier into the middle of the box and move it to the right at the velocity $v$ for a unit time, or (iii) do nothing. At the end of control actions, the barrier is removed from the engine. We assume that the insertion and removal of the barrier  is frictionless and as such do not consume any work. The sequence of operations is depicted in Fig.~\ref{fig:szilard}.
Denote by $p(x_{t+1}|x_t,u_t)$ the transition probability from the state $x_t$ to another state $x_{t+1}$ when control action $u_t$ is applied.
By $\mathbb{E}c(X_t,U_t)$ we denote the expected work required to apply control action $u_t$ at time $t$ when the state of the engine is $x_t$.\footnote{Here, we do not provide a detailed model of the function $c(x_t,u_t)$. See for instance \cite{horowitz2014second} for a model of work extraction based on the Langevin equation.} This quantity is negative if the controller is expected to extract work from the engine. Work extraction occurs when the gas molecule collides with the barrier and ``pushes'' it in the direction of its movement.

Right before applying a control action $U_t$, suppose that the controller makes (a possibly noisy) observation of the engine state, and thus there is an information flow from $X^t$ to $U^t$. For our discussion, there is no need to describe what kind of sensing mechanism is involved in this step. However, notice that if an error-free observation of the engine state $X_t$ is performed, then the controller can choose a control action such that $\mathbb{E}c(X_t,U_t)$ is always non-positive. (Consider moving the barrier always to the opposite direction from the position of the gas molecule.)
At first glance, this seems to imply that one can construct a device that is expected to cyclically extract work from a single thermal bath, which is a contradiction to the Kelvin-Planck statement of the second law of thermodynamics.

\ifdefined\IFONECOL
\begin{figure}[t]
    \centering
    \includegraphics[width=0.6\columnwidth]{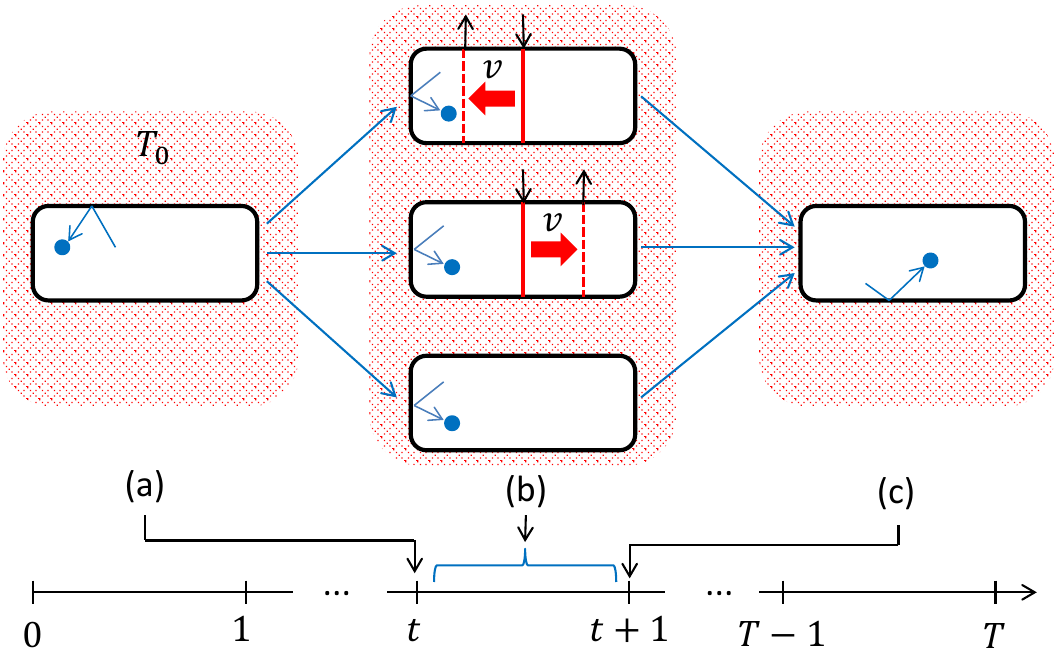}
    \caption{Modified Szilard engine. The controller performs the following steps in a unit time. (a) The controller makes (a possibly noisy) observation of the state $X_t$ of the engine. (b) One of the three possible control actions $U_t$ is applied. (c) Barrier is removed.}
    \label{fig:szilard}
\end{figure}
\else
\begin{figure}[t]
    \centering
    \includegraphics[width=\columnwidth]{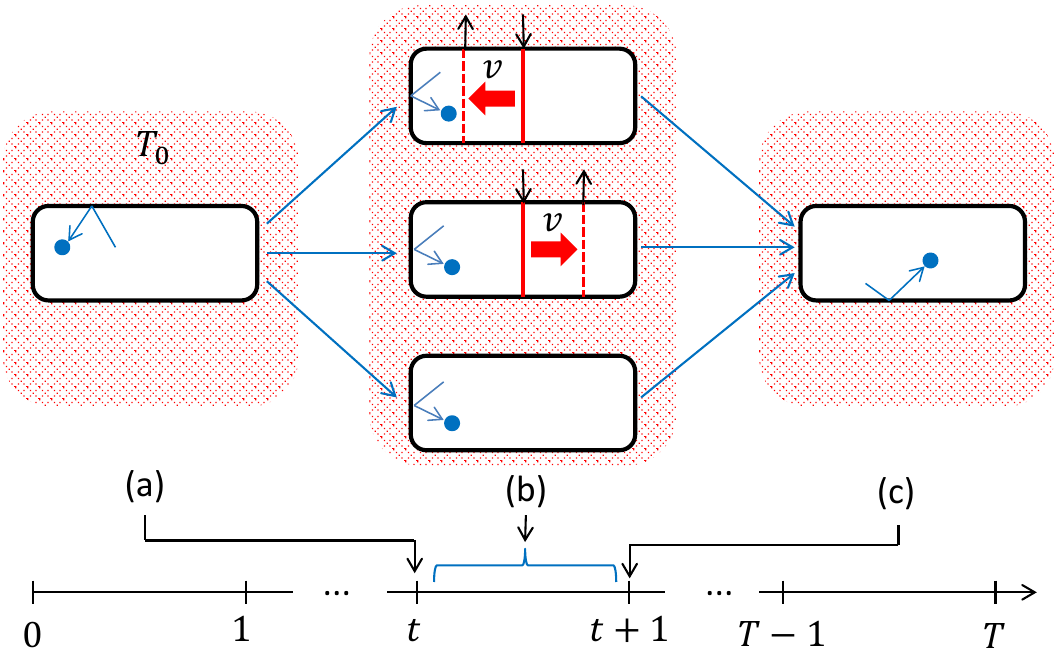}
    \caption{Modified Szilard engine. The controller performs the following steps in a unit time. (a) The controller makes (a possibly noisy) observation of the state $X_t$ of the engine. (b) One of the three possible control actions $U_t$ (move the barrier to the left or to the right, or do nothing) is applied. (c) At the end of control action, the barrier is removed.}
    \label{fig:szilard}
\end{figure}
\fi

\begin{figure*}[t]
\begin{minipage}{0.35\textwidth}
    \centering
    \includegraphics[width=0.65\textwidth]{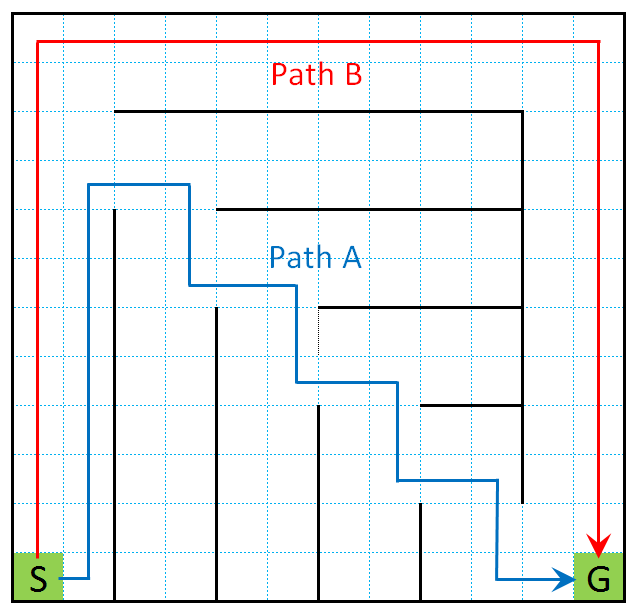}
    \caption{Information-regularized optimal navigation through a maze.}
    \label{fig:grid}
    \vspace{2ex}
    \includegraphics[width=0.9\textwidth]{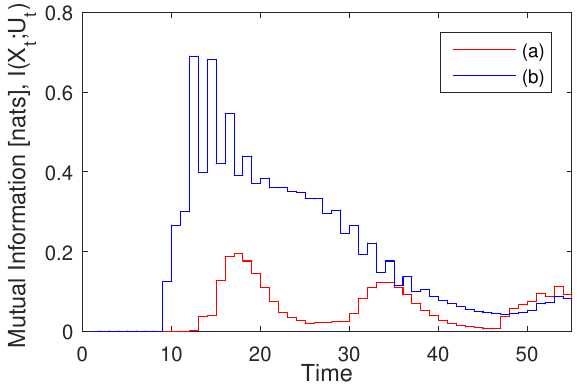}
    \vspace{-1ex}
    \caption{Information usage by the policies (a) and (b) in Fig.~\ref{fig:maze}.}
    \label{fig:mi}
\end{minipage}
\begin{minipage}{0.63\textwidth}
    \centering
    \includegraphics[width=\textwidth]{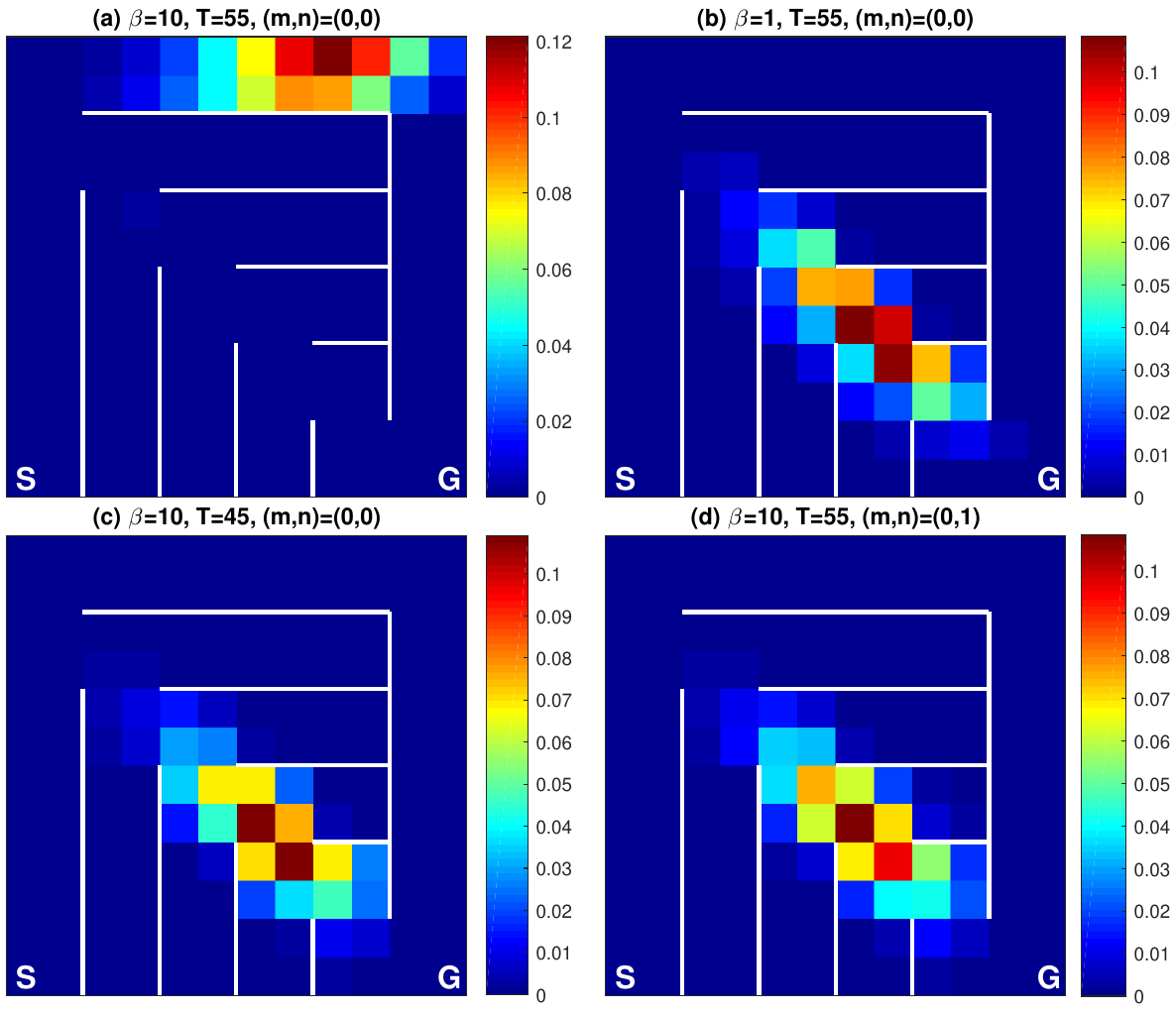}
    \caption{State probability distribution $\mu_t(x_t)$ at $t=25$.}
    \label{fig:maze}
\end{minipage}
\end{figure*}

It is now widely recognized that this paradox (Maxwell's demon) can be resolved  by including the ``memory'' of the controller into the picture.
Recently, a generalized second law is proposed by \cite{ito2013information}, in which transfer entropy plays a critical role.
Viewing the combined engine and memory system as a Bayesian network comprised of $X_t$ and $U_t$ (see \cite{ito2013information} for details), and
assuming that the free energy change of the engine from $t=0$ to $t=T$ is zero (which is the case when the above sequence of operations are repeated in a cyclic manner with period $T$), the generalized second law \cite[equation (10)]{ito2013information} reads
\begin{equation}
\label{eqsecondlaw}
\sum_{t=0}^{T-1} \mathbb{E}c(X_t,U_t)+k_BT_0 I(X_0^{T-1}\rightarrow U_0^{T-1}) \geq 0
\end{equation}
where $k_B$ [J/K] is the Boltzmann constant.
The above inequality shows that a positive amount of work is extractable (i.e., the first term can be negative), but this is possible only at the expense of the transfer entropy cost (the second term must be positive).\footnote{The consistency with the classical second law is maintained if one accepts Landauer's principle, which asserts that erasure of one bit of information from any sort of memory device in an environment at temperature $T_0$ [K] requires at least $k_B T_0 \log 2$ [J] of work. See \cite{bennett1982thermodynamics} for the further discussions.}
Given a fundamental law \eqref{eqsecondlaw}, a natural question is how efficient the considered thermal engine can be by optimally designing a control policy $q(u_t|x^t,u^{t-1})$. This can be analyzed by minimizing the left hand side of \eqref{eqsecondlaw}, which is precisely the TERMDP problem \eqref{eqmainprob}.

\section{Numerical experiment}
\label{secnum}

In this section, we apply the proposed forward-backward Arimoto-Blahut algorithm (Algorithm~\ref{alg1}) to study how the price of information affects the level of information-frugality, which yields qualitatively different decision policies.

Consider a situation in which Alice, whose movements are described by Markovian dynamics controlled by Bob, is traveling through a maze shown in Fig.~\ref{fig:grid}.
Suppose that Bob knows the geometry of the maze (including start and goal locations), but observing Alice's location is costly.
We model this problem as an MDP where the state $X_t$ is the cell where Alice is located at time step $t$, and $U_t$ is a navigation instruction given by Bob. The observation cost is characterized by the transfer entropy. We assume five different instructions are possible; $u=N, E, S, W$ and $R$, corresponding to \emph{go north}, \emph{go east}, \emph{go south}, \emph{go west}, and \emph{rest}. The initial state is the cell indicated by ``S'' in Fig.~\ref{fig:grid}, and the motion of Alice is described by a transition probability $p(x_{t+1}|x_t, u_t)$.

The transition probability is defined by the following rules.
At each cell, a transition to the indicated direction occurs w.p. $0.8$ if there is no wall in the indicated direction, while transitions to any open directions (directions without walls) occurs w.p. $0.05$ each. With the remaining probability, Alice stays in the same cell. If there is a wall in the indicated direction, or $u=R$, then transition to each open direction occurs w.p. $0.05$, while Alice stays in the same cell with the remaining probability.

At each time step $t=1, 2, ... , T$, the state-dependent cost is defined by $c_t(x_t,u_t)=0$ if $x_t$ is already the target cell indicated by ``G'' in Fig.~\ref{fig:grid}, and $c_t(x_t,u_t)=1$ otherwise.
The terminal cost is $0$ if $x_{T+1}=G$ and $10000$ otherwise.
We consider transfer entropy $I_{m,n}(X^T\rightarrow U^T)$ to quantify the price of information that Bob must acquire about Alice's location. With some nonnegative weight $\beta$, the overall control problem can be written as  \eqref{eqmainprob}.

As shown in Fig.~\ref{fig:grid}, there are two qualitatively different paths from the origin to the target. The path A is shorter than the path B, and hence Bob will try to navigate Alice along path A when no information-theoretic cost is considered (i.e., $\beta=0$). However, navigating Alice along the path A requires Bob to have accurate information about Alice's current location, as this path is ``risky'' (there are many side roads with dead ends).
The path B is longer, but navigating through it is simpler; rough knowledge about Alice's location is sufficient to provide correct instructions.
Hence, it is expected that Bob would try to navigate Alice through A when information is relatively cheap ($\beta$ is small), while he would choose B when information is expensive ($\beta$ is large).

Fig.~\ref{fig:maze} shows the solutions to the considered problem. Solutions are obtained by iterating Algorithm~\ref{alg1} sufficiently many times in four different conditions. Each plot shows a snapshot of the state probability distribution $\mu_t(x_t)$ at time $t=25$. Fig.~\ref{fig:maze} (a) is obtained under the setting that the cost of information is high ($\beta=10$), the planning horizon is long ($T=55$), and the transfer entropy of degree $(m,n)=(0,0)$ is considered. Accordingly, a decision policy of the form of $q_t(u_t|x_t)$ is considered. It can be seen that with high probability, the agent is navigated through the longer path. In Fig.~\ref{fig:maze} (b), the cost of information is reduced  ($\beta=1$) while the other settings are kept the same.
As expected, the solution chooses the shorter path.
Fig.~\ref{fig:mi} shows the time-dependent information usage in (a) and (b); it shows that the total information usage is greater in situation (b) than in (a).

We note that this simulation result is consistent with a prior work \cite{rubin2012trading}, where similar numerical experiments were conducted.
Using Algorithm~\ref{alg1}, we can further investigate the nature of the problem. Fig.~\ref{fig:maze} (c) considers the same setting as in (a) except that the planning horizon is shorter ($T=45$). This result shows that the solution becomes qualitatively different depending on how close the deadline is even if the cost of information is the same. Finally, Fig.~\ref{fig:maze} (d) considers the case where the transfer entropy has degree $(m,n)=(0,1)$ and the decision policy is of the form of $q_t(u_t|x_t,u_{t-1})$. Although the rest of simulation parameters are unchanged from (a), we observe that the shorter path is chosen in this case.
This result demonstrates that the solution to \eqref{eqmainprob} can be qualitatively different depending on the considered degree of transfer entropy costs.

\section{Summary and Future Work}
\label{secsummary}

{\updnew
In this paper, we considered a mathematical framework of transfer-entropy-regularized Markov Decision Process (TERMDP), which is motivated both in engineering (networked
  control systems) and scientific (non-equilibrium thermodynamics) contexts.  
We derived structural properties of the optimal solution, and provided a necessary optimality condition written  as a set of coupled nonlinear equations.
We proposed an iterative numerical algorithm (forward-backward Arimoto-Blahut algorithm) in which every limit point of the generated sequence is guaranteed to be a stationary point of the given TERMDP. By numerical simulation (information-constrained maze navigation), it was demonstrated that the proposed algorithm can be used to find an optimal solution candidate.

The proposed algorithm has several limitations as summarized in Remark~\ref{rem1}, which must be addressed in future work. 
Improvement of the convergence speed of the proposed algorithm is also necessary for many applications.
Finally, the roles of transfer entropy in engineering and scientific applications (including networked control systems, non-equilibrium thermodynamics and beyond) and implications of TERMDP solutions studied in this paper need further investigation.
}

%

{\upd
\section*{Acknowledgement}                            
The authors would like to thank the anonymous reviewers for their valuable suggestions.
}

\ifdefined\LONGVERSION

\appendix

\subsection{Proof of Proposition~\ref{propstructure}}
\label{appstructure}
Our proof is based on backward induction.

\underline{$k=T$}: We prove (a) first. For a given policy $q_T(u_t|x^T, u^{T-1})$, let 
$
\lambda_T(x^T, u^T)\triangleq q_T(u_T|x^T, u^{T-1})\mu_T(x^T, u^{T-1})
$
be the joint distribution induced by $q_T$.
{\updnew
Notice that 
\begin{equation}
\label{eqlambdamu}
\lambda_T(x^T, u^{T-1})=\mu_T(x^T,u^{T-1})
\end{equation}
holds by construction.
}
Let $\lambda_T(x_T, u_{T-n}^T)$ and $\lambda_T(x_T, u_{T-n}^{T-1})$ be marginals of $\lambda_T(x^T, u^T)$. Construct a new policy $q'_T$ as
{\updnew
\begin{equation}
\label{eqnewqconstruct}
q'_T(u_T|x_T, u_{T-n}^{T-1})\triangleq \frac{\lambda_T(x_T, u_{T-n}^T)}{\lambda_T(x_T, u_{T-n}^{T-1})}
\end{equation}
if $\lambda_T(x_T, u_{T-n}^{T-1})>0$, and as an arbitrary probability distribution on $\mathcal{U_T}$ if $\lambda_T(x_T, u_{T-n}^{T-1})=0$.
}
Let 
\begin{equation}
\label{eqlambdapdef}
\lambda'_T(x^T, u^T)\triangleq q'_T(u_T|x_T, u_{T-n}^{T-1})\mu_T(x^T, u^{T-1})
\end{equation}
be the joint distribution induced by $q'_T$. Then, we have
\begin{equation}
\label{eqlambdap}
\lambda_T(x_T, u_{T-n}^T)=\lambda'_T(x_T, u_{T-n}^T), 
\end{equation}
which can be directly verified as
\begin{subequations}
\begin{align}
&\lambda'_T(x_T, u_{T-n}^T)\nonumber \\
&=\sum_{\substack{ x^{T-1}\in\mathcal{X}^{T-1} \\ u^{T-n-1}\in\mathcal{U}^{T-n-1} }} q'_T(u_T|x_T, u_{T-n}^{T-1})\mu_T(x^T, u^{T-1}) \label{eqlambda1}\\
&=\sum_{\substack{ x^{T-1}\in\mathcal{X}^{T-1} \\ u^{T-n-1}\in\mathcal{U}^{T-n-1} }} q'_T(u_T|x_T, u_{T-n}^{T-1})\lambda_T(x^T, u^{T-1}) \label{eqlambda2}\\
&=q'_T(u_T|x_T, u_{T-n}^{T-1})\lambda_T(x_T, u_{T-n}^{T-1}) \nonumber \\
&=\lambda_T(x_T, u_{T-n}^T) \label{eqlambda4}
\end{align}
\end{subequations}
{\updnew where the equality \eqref{eqlambda1} holds by definition \eqref{eqlambdapdef},  \eqref{eqlambda2} by \eqref{eqlambdamu},  and \eqref{eqlambda4} by the construction \eqref{eqnewqconstruct}.}
Now the Bellman equation at $k=T$ reads
\[
V_T(\mu_T(x^T,u^{T-1}))=\min_{q_T} J_T^c(\lambda_T)+J_T^I(\lambda_T)
\]
where
\begin{align*}
J_T^c(\lambda_T)&=\mathbb{E}^{\lambda_T, p_{T+1}}\left(c_T({\upd X_T, U_T})+c_{T+1}({\upd  X_{T+1}})\right) \\
J_T^I(\lambda_T)&=I_{\lambda_T}(X_{T-m}^T;U_T|U_{T-n}^{T-1}).
\end{align*}
To establish (a) for $k=T$, it is sufficient to show that $J_T^c(\lambda_T)=J_T^c(\lambda'_T)$ and $J_T^I(\lambda_T)\geq J_T^I(\lambda'_T)$. The first equality holds because of \eqref{eqlambdap}. To see the second inequality,
\begin{subequations}
\begin{align}
J_T^I(\lambda_T)&=I_{\lambda_T}(X_{T-m}^T;U_T|U_{T-n}^{T-1})  \\
&\geq  I_{\lambda_T}(X_T;U_T|U_{T-n}^{T-1}) \\
&= I_{\lambda'_T}(X_T;U_T|U_{T-n}^{T-1}) \label{eqJI1}\\
&= I_{\lambda'_T}(X_T;U_T|U_{T-n}^{T-1}) \nonumber \\
&\hspace{10ex}+I_{\lambda'_T}(X_{T-m}^{T-1};U_T|X_T, U_{T-n}^{T-1}) \label{eqJI2} \\
&=I_{\lambda'_T}(X_{T-m}^T;U_T|U_{T-n}^{T-1}) \\
&=J_T^I(\lambda'_T).
\end{align}
\end{subequations}
The equality \eqref{eqJI1} follows from \eqref{eqlambdap}.
The second term in \eqref{eqJI2} is zero since $U_T$ is independent of $X_{T-m}^{T-1}$ given $(X_T, U_{T-n}^{T-1})$ by construction of $\lambda'_T$ in \eqref{eqlambdapdef}. Hence the statement (a) is established for $k=T$.

The statement (b) follows immediately from the above discussion. To establish (c), notice that due to (a), we can assume the optimal policy of the form $q_T(u_T|x_T, u_{T-n}^{T-1})$ without loss of generality. Hence, due to (b), the Bellman equation at $k=T$ can be written as
\begin{align*}
V_T(\mu_T(x^T, u^{T-1}))=&\min_{q_T(u_T|x_T, u_{T-n}^{T-1})} \Bigl\{I(X_T;U_T|U_{T-n}^{T-1})\Bigr. \\
&+\Bigl.\mathbb{E} c_T({\upd X_T, U_T})+\mathbb{E}  c_{T+1}({\upd X_{T+1}}) \Bigr\}.
\end{align*}
Since the last expression depends on $\mu_T(x^T, u^{T-1})$ only through its marginal $\mu_T(x_T, u_{T-n}^{T-1})$, we conclude that $V_T(\mu_T(x^T, u^{T-1}))$ is a function of the marginal $\mu_T(x_T, u_{T-n}^{T-1})$ only. This establishes (c) for $k=T$.

\underline{$k=t$}: Next, we assume (a), (b) and (c) hold for $k=t+1$. To establish (a), (b) and (c) for $k=t$, notice that under the induction hypothesis, we can assume without loss of generality that policies for $k=t+1, t+2, ... , T$ are of the form $q_k(u_k|x_k,u_{k-n}^{k-1})$, and that the value function at $k=t+1$ depends only on $\mu_{t+1}(x_{t+1}, u_{t-n+1}^t)$. With a slight abuse of notation, the latter fact is written as
\[
V_{t+1}(\mu_{t+1}(x^{t+1}, u^t))=V_{t+1}(\mu_{t+1}(x_{t+1}, u_{t-n+1}^t)).
\]
Thus, the Bellman equation at $k=t$ can be written as
\begin{align}
&V_t\left(\mu_t(x^t, u^{t-1})\right) 
=\min_{q_t(u_t|x^t,u^{t-1})} \Bigl\{ \mathbb{E}c_t({\upd X_t, U_t}) \biggr.  \nonumber \\
& \hspace{1ex}\Bigl. + I(X_{t-m}^t;U_t|U_{t-n}^{t-1})  + V_{t+1}(\mu_{t+1}(x_{t+1}, u_{t-n+1}^{t}))\Bigr\}. \label{eqvfindt}
\end{align}
Now, using the similar construction to the case for $k=T$, one can show that for every $q_t(u_t|x^t,u^{t-1})$, there exists a policy of the form $q'_t(u_t|x_t,u_{t-n}^{t-1})$ such that the value of the objective function in the right hand side of \eqref{eqvfindt} attained by $q'_t$ is less than or equal to the value attained by $q_t$. This observation establishes (a) for $k=t$. Statements (b) and (c) for $k=t$ follows similarly.

{\updnew
\ifdefined\IFONECOL
\section{Local minima of simplified TERMDP}
\else
\subsection{Local minima of simplified TERMDP}
\fi
\label{apdxlocal}
Let $q=\{q(u_t|x^t,u^{t-1})\}_{t=1}^T$ be a given policy.
We say that the joint distribution $\lambda(x^T, u^T)$ is \emph{induced} by $q$ if it is defined as
\[
\lambda(x^T,u^T)=\prod_{t=1}^T q(u_t|x^t, u^{t-1})p(x_t|x_{t-1},u_{t-1})
\]
with the initial time condition $p(x_1|x_0, u_0)=\mu(x_1)$.
As in the proof of Proposition~\ref{propstructure}, we construct a new policy $q'=\{q'(u_t|x_t,u_{t-n}^{t-1})\}_{t=1}^T$ from $\lambda(x^T, u^T)$ by
\[
q'(u_t|x_t, u_{t-n}^{t-1})=\frac{\lambda(x_t, u_{t-n}^t)}{\lambda(x_t, u_{t-n}^{t-1})}
\]
for each $t=1, 2, ... , T$. For simplicity of the analysis, we assume that $\lambda(x^T, u^T)$ satisfies the following condition:
\begin{assumption}
\label{asmp1}
For each $t=1, 2,  ... , T$ and $(x_t, u_{t-n}^{t-1})\in\mathcal{X}_t\times \mathcal{U}_{t-n}^{t-1}$, we have $\lambda(x_t, u_{t-n}^{t-1})>0$.
\end{assumption}
Although this assumption is somewhat restrictive, it is valid when, for instance, the transition probability $p(x_t|x_{t-1},u_{t-1})$ and the policy $q(u_t|x^t, u^{t-1})$ assign a nonzero probability mass to each element of $\mathcal{X}_t$ and $\mathcal{U}_t$, respectively.
Denote by $\lambda'(x^T, u^T)$ the joint distribution induced by $q'$, i.e.,
\begin{align}
\lambda'(x^T,u^T)&=\prod_{t=1}^T q'(u_t|x^t, u^{t-1})p(x_t|x_{t-1},u_{t-1}) \nonumber \\
&=\prod_{t=1}^T \frac{\lambda(x_t, u_{t-n}^t)}{\lambda(x_t, u_{t-n}^{t-1})}p(x_t|x_{t-1},u_{t-1}). \label{eqlambdarational}
\end{align}
We write the construction of $q'$ from $q$ via the procedure above as $q'=\pi(q)$, where $\pi:\mathcal{Q}\rightarrow\mathcal{Q}'$ is a mapping from the space $\mathcal{Q}$ of general policies of the form $\{q(u_t|x^t,u^{t-1})\}_{t=1}^T$ to the space $\mathcal{Q}'$ of simplified policies of the form $\{q(u_t|x_t,u_{t-n}^{t-1})\}_{t=1}^T$. Notice that $\mathcal{Q}'\subset \mathcal{Q}$. 

We first show that $\pi$ is a continuous mapping with respect to an appropriate metric on $\mathcal{Q}$. Suppose $\bar{q}\in\mathcal{Q}$ and $\tilde{q}\in\mathcal{Q}$ are given policies, and let $\bar{\lambda}$ and  $\tilde{\lambda}$ be joint distributions induced by $\bar{q}$ and $\tilde{q}$, respectively.
We define a metric $\|\bar{q}-\tilde{q}\|$ as the $\ell_1$ distance between $\bar{\lambda}$ and $\tilde{\lambda}$:
\begin{equation}
\label{eqqmetric}
\|\bar{q}-\tilde{q}\|\triangleq \sum_{x^T\in\mathcal{X}^T, u^T\in \mathcal{U}^T} |\bar{\lambda}(x^T, u^T)-\tilde{\lambda}(x^T, u^T)|.
\end{equation}
\begin{claim}
\label{claim1}
Suppose $\bar{q}\in \mathcal{Q}$ and $\tilde{q}\in \mathcal{Q}$ satisfy Assumption~\ref{asmp1}. For every $\epsilon>0$, there exists $\delta>0$ such that
\[
\|\bar{q}-\tilde{q}\|<\delta \Rightarrow \|\pi(\bar{q})-\pi(\tilde{q})\|<\epsilon.
\]
\end{claim}
\begin{pf}
Let $\bar{\lambda}'$ and $\tilde{\lambda}'$ be joint distributions induced by $\bar{q}'=\pi(\bar{q})$ and $\tilde{q}'=\pi(\tilde{q})$, respectively. Notice that $\|\bar{q}-\tilde{q}\|<\delta$ means that $\sum_{\mathcal{X}^T, \mathcal{U}^T} |\bar{\lambda}(x^T, u^T)-\tilde{\lambda}(x^T, u^T)|<\delta$, which implies that
\begin{align*}
&|\bar{\lambda}(x_t, u_{t-n}^{t})-\tilde{\lambda}(x_t, u_{t-n}^{t})|<\delta \text{ and } \\
&|\bar{\lambda}(x_t, u_{t-n}^{t-1})-\tilde{\lambda}(x_t, u_{t-n}^{t-1})|<\delta.
\end{align*}
Since \eqref{eqlambdarational} shows that $\lambda'(x^T, u^T)$ for each $(x^T, u^T)\in\mathcal{X}^T\times \mathcal{U}^T$ is a rational function of $\lambda(x_t, u_{t-n}^{t})$ and $\lambda(x_t, u_{t-n}^{t-1})$, and a rational function is continuous in the domain of strictly positive denominators, for each $\epsilon>0$ and $(x^T, u^T)\in\mathcal{X}^T\times \mathcal{U}^T$, there exists $\delta>0$ such that $\|\bar{q}-\tilde{q}\|<\delta$ implies 
\[
|\bar{\lambda}'(x^T, u^T)-\tilde{\lambda}'(x^T, u^T)|<\frac{\epsilon}{|\mathcal{X}^T||\mathcal{U}^T|},
\]
which further implies 
\begin{align*}
\|\pi(\bar{q})-\pi(\tilde{q})\|
&=\|\bar{q}'-\tilde{q}'\| \\
&=\sum_{x^T\in \mathcal{X}^T, u^T\in \mathcal{U}^T}  |\bar{\lambda}'(x^T, u^T)-\tilde{\lambda}'(x^T, u^T)| \\
&< \sum_{x^T\in \mathcal{X}^T, u^T\in \mathcal{U}^T} \frac{\epsilon}{|\mathcal{X}^T||\mathcal{U}^T|} \\
&=\epsilon.
\end{align*}
\end{pf}
The next claim shows that $\pi$ leaves an element of $\mathcal{Q}'$ invariant.
\begin{claim}
\label{claim2}
Suppose $q\in\mathcal{Q}'$ and the induced joint distribution $\lambda$ satisfies Assumption~\ref{asmp1}. Then $\pi(q)=q$.
\end{claim}
\begin{pf}
By definition, 
\[
\lambda(x^t, u^t)=\prod_{k=1}^t q(u_k|x_k, u_{k-n}^{k-1})p(x_k|x_{k-1},u_{k-1}).
\]
Therefore,
\begin{align}
&\lambda(x_t, u_{t-n}^t)=\!\!\!\!\!\!\!\!\sum_{\substack{x^{t-1}\in\mathcal{X}^{t-1} \\ u^{t-n-1} \mathcal{U}^{t-n-1}}}\prod_{k=1}^t q(u_k|x_k, u_{k-n}^{k-1})p(x_k|x_{k-1},u_{k-1}) \nonumber \\
&=q(u_t|x_t, u_{t-n}^{t-1})\sum_{x^{t-1}\in\mathcal{X}^{t-1}}p(x_t|x_{t-1},u_{t-1}) \nonumber \\
&\hspace{4ex} \times \sum_{u^{t-n-1}\in\mathcal{U}^{t-n-1}} \prod_{k=1}^{t-1} q(u_k|x_k, u_{k-n}^{k-1})p(x_k|x_{k-1},u_{k-1}) \label{eqlambdamarginal1}
\end{align}
and
\begin{align}
&\lambda(x_t, u_{t-n}^{t-1})=\sum_{u_t \in \mathcal{U}_t}\lambda(x_t, u_{t-n}^t) \nonumber \\
&=\sum_{x^{t-1}\in\mathcal{X}^{t-1}}p(x_t|x_{t-1},u_{t-1}) \nonumber \\
&\hspace{4ex} \times \sum_{u^{t-n-1}\in\mathcal{U}^{t-n-1}} \prod_{k=1}^{t-1} q(u_k|x_k, u_{k-n}^{k-1})p(x_k|x_{k-1},u_{k-1}).  \label{eqlambdamarginal2}
\end{align}
Let $q'=\pi(q)$. Then, from \eqref{eqlambdamarginal1} and  \eqref{eqlambdamarginal2}, we have
\[
q'(u_t|x_t, u_{t-n}^{t-1})=\frac{\lambda(x_t, u_{t-n}^t)}{\lambda(x_t, u_{t-n}^{t-1})} = q(u_t|x_t, u_{t-n}^{t-1}).
\]
Therefore, $\pi(q)=q$.
\end{pf}

We are now ready to prove the following proposition stating that a local minimum for the simplified TERMDP \eqref{eqsimpleprob} is a local minimum for the original TERMDP \eqref{eqmainprob}. For a given policy $q\in\mathcal{Q}$, denote by
\[
f(q)=J(X^{T+1},U^T)+\beta I_{m,n}(X^T\rightarrow U^T)
\]
the value of the TERMDP \eqref{eqmainprob}. By Proposition~\ref{propstructure}, it also follows that if $q\in\mathcal{Q}'$, then $f(q)$ is equal to the value of the simplified TERMDP \eqref{eqsimpleprob}, i.e., 
\[
f(q)=J(X^{T+1},U^T)+\beta I_{0,n}(X^T\rightarrow U^T).
\]
\begin{prop}
Suppose $q^*=\{q^*(u_t|x_t, u_{t-n}^{t-1})\}_{t=1}^T\in\mathcal{Q}'$ satisfies Assumption~\ref{asmp1}. If $q^*$ is a local minimizer for the simplified TERMDP \eqref{eqsimpleprob}, i.e.,
\begin{itemize}
\item[(a)] there exists $\epsilon>0$ such that $f(q^*)\leq f(q')$ holds for all $q'\in\mathcal{Q}'$ with $\|q'-q^*\|<\epsilon$.
\end{itemize} 
Then, $q^*$ is a local minimizer for the original TERMDP \eqref{eqmainprob}, i.e.,
\begin{itemize}
\item[(b)] there exists $\delta>0$ such that $f(q^*)\leq f(q)$ holds for all $q\in\mathcal{Q}$ with $\|q-q^*\|<\delta$.
\end{itemize} 
\end{prop}
\begin{pf}
Pick a constant $\epsilon>0$ such that the condition (a) holds. 
By continuity of $\pi$ (Claim~\ref{claim1}), there exists a constant $\delta>0$ such that
\begin{equation}
q\in\mathcal{Q}, \|q-q^*\|<\delta \Rightarrow \|\pi(q)-q^*\|<\epsilon. \label{eqprop6proof1}
\end{equation}
Here, we have used the fact that $\pi(q^*)=q^*$ (Claim~\ref{claim2}).

To complete the proof by contradiction, suppose the negation of the condition (b) holds:
\begin{itemize}
\item[($\neg$b)] for every $\delta>0$, there exists $\bar{q}\in\mathcal{Q}$ such that $\|\bar{q}-q^*\|<\delta$ and $f(\bar{q})<f(q^*)$.
\end{itemize} 
Now, pick a policy $\bar{q}\in\mathcal{Q}$ such that $\|\bar{q}-q^*\|<\delta$ and 
\begin{equation}
 \label{eqprop6proof2}
f(\bar{q})<f(q^*).
\end{equation}
If we write $\bar{q}'\triangleq \pi(\bar{q})\in \mathcal{Q}'$, it also follows from \eqref{eqprop6proof1} that 
\begin{equation}
 \label{eqprop6proof3}
\|\bar{q}'-q^*\|<\epsilon.
\end{equation}
By Proposition~\ref{propstructure},  \eqref{eqprop6proof2} implies that
\begin{equation}
\label{eqprop6proof4}
f(\bar{q}')\leq f(\bar{q})<f(q^*).
\end{equation}
However, \eqref{eqprop6proof3} and \eqref{eqprop6proof4} contradict (a). Therefore, the condition (b) must hold.
\end{pf}
}

\subsection{Proof of Proposition~\ref{propfinitete}}
\label{apptemonotone}
For each $n'\geq n$,
\begin{align*}
&I_{0,n'}(X^T\rightarrow U^T) =\sum\nolimits_{t=1}^T I(X_t;U_t|U_{t-n'}^{t-1}) \\
&=\sum\nolimits_{t=1}^T H(U_t|U_{t-n'}^{t-1})-H(U_t|X_t,U_{t-n'}^t) \\
&=\sum\nolimits_{t=1}^T H(U_t|U_{t-n'}^{t-1})-H(U_t|X_t,U_{t-n}^t).
\end{align*}
In the last step, we used the fact that $H(X|Y,Z)=H(X|Y)$ holds when $X$ and $Z$ are conditionally independent given $Y$. By the structure of $q_t(u_t|x_t, u_{t-n}^{t-1})$, $U_t$ and $(X_t, U_{t-n}^{t-1})$ is conditionally independent of $U_{t-n'}^{t-n-1}$. Now,
\begin{align*}
&I_{0,n'}(X^T\rightarrow U^T)-I_{0,n'+1}(X^T\rightarrow U^T) \\
&=\sum\nolimits_{t=1}^T \left(H(U_t|U_{t-n'}^{t-1})-H(U_t|U_{t-n'-1}^{t-1})\right) \geq 0
\end{align*}
since entropy never increases by conditioning.

\ifdefined\IFONECOL
\section{Proof of Proposition~\ref{propcom2}}
\else
\subsection{Proof of Proposition~\ref{propcom2}}
\fi
\label{app1}

For each $t=1,2, ... , T$, we have
\ifdefined\IFONECOL
\begin{align*}
I(X_{t-m}^t; W_t|U_{t-n}^{t-1}) &=I(X_{t-m}^t; W_t, U_t|U_{t-n}^{t-1}) \\
&=I(X_{t-m}^t; U_t|U_{t-n}^{t-1})+I(X_{t-m}^t;W_t|U_{t-n}^t) \\
&\geq I(X_{t-m}^t; U_t|U_{t-m}^{t-1}).
\end{align*}
\else
\begin{align*}
&I(X_{t-m}^t; W_t|U_{t-n}^{t-1}) \\
&=I(X_{t-m}^t; W_t, U_t|U_{t-n}^{t-1}) \\
&=I(X_{t-m}^t; U_t|U_{t-n}^{t-1})+I(X_{t-m}^t;W_t|U_{t-n}^t) \\
&\geq I(X_{t-m}^t; U_t|U_{t-m}^{t-1}).
\end{align*}
\fi
The first equality is due to the particular structure of the decoder specified by (A) or (B). Thus
\[
\sum\nolimits_{t=1}^T I(X_{t-m}^t; W_t|U_{t-n}^{t-1}) \geq I_{m,n}(X^T\rightarrow U^T).
\]
The proof of Proposition~\ref{propcom2} is complete by noticing the following chain of inequalities.
\begin{align*}
R\log 2\! &=\!\sum\nolimits_{t=1}^T R_t \log 2 \\
&\geq \!\sum\nolimits_{t=1}^T \!H(W_t) \\
&\geq \!\sum\nolimits_{t=1}^T \!H(W_t|W^{t-1},U_{t-n}^{t-1}) \\
&\geq \!\sum\nolimits_{t=1}^T \!H(W_t|U_{t-n}^{t-1})-H(W_t|X_{t-m}^t, U_{t-n}^{t-1}) \\
&=\!\sum\nolimits_{t=1}^T \!I(X_{t-m}^t;W_t|U_{t-n}^{t-1}).
\end{align*}

\fi







\bibliographystyle{IEEEtran}
\bibliography{refs}

\end{document}